\documentclass{amsart}

\usepackage{mystyle}

\makeatletter
\def\@tocline#1#2#3#4#5#6#7{\relax
  \ifnum #1>\c@tocdepth 
  \else
    \par \addpenalty\@secpenalty\addvspace{#2}%
    \begingroup \hyphenpenalty\@M
    \@ifempty{#4}{%
      \@tempdima\csname r@tocindent\number#1\endcsname\relax
    }{%
      \@tempdima#4\relax
    }%
    \parindent\z@ \leftskip#3\relax \advance\leftskip\@tempdima\relax
    \rightskip\@pnumwidth plus4em \parfillskip-\@pnumwidth
    #5\leavevmode\hskip-\@tempdima
      \ifcase #1
       \or\or \hskip 1em \or \hskip 2em \else \hskip 3em \fi%
      #6\nobreak\relax
    \dotfill\hbox to\@pnumwidth{\@tocpagenum{#7}}\par
    \nobreak
    \endgroup
  \fi}
\makeatother

\begin{document}

\title{Exotic Motivic Periodicities}

\author{Bogdan Gheorghe}
\email{gheorghebg@wayne.edu}

\address{Department of Mathematics\\ Wayne State University\\
Detroit, MI 48202, USA}
\date{\today}


\begin{abstract}
One can attempt to study motivic homotopy groups by mimicking the classical (non-motivic) chromatic approach. There are however major differences, which makes the motivic story more complicated and still not well understood. For example, classically the $p$-local sphere spectrum $S^0_{(p)}	$ admits an essentially unique non-nilpotent self-map, which is not the case motivically, since Morel showed that the first Hopf map $\eta \colon S^{1,1} \lto \Ss$ is non-nilpotent. In the same way that the non-nilpotent self-map $2 = v_0 \in \piss(\Ss)$ starts the usual chromatic story of $v_n$-periodicity, there is a similar theory starting with the non-nilpotent element $\eta \in \piss(\Ss)$, which Andrews-Miller denoted by $\eta = w_0$. In this paper we investigate the beginning of the motivic story of $w_n$-periodicity when the base scheme is $\mathbf{Spec} \! \ \mathbb{C}$. In particular, we construct motivic fields $K(w_n)$ designed to detect such $w_n$-periodic phenomena, in the same way that $K(n)$ detects $v_n$-periodic phenomena. In the hope of detecting motivic nilpotence, we also construct a more global motivic spectrum $\wBP$ with homotopy groups $\piss(\wBP) \cong \F_2[w_0, w_1, \ldots]$.

\end{abstract}

\maketitle

\tableofcontents

\section{Introduction}

\subsection{Motivation}

The \emph{chromatic approach} to classical homotopy theory is a very powerful organizational tool to study the homotopy category of ($p$-local) finite CW-complexes. In particular, given a $p$-local finite CW-complex $X$, the chromatic approach provides an algorithm for computing its homotopy groups $\pi_{\ast}(X)$. This algorithm relies heavily on the existence of the \emph{chromatic filtration} on $\pi_{\ast}(X)$. This is an increasing filtration indexed by non-negative integers $n \in \mathbb{N}_{0}$, where we say that an element $f \in \pi_{\ast}(X)$ in the $n^{\text{th}}$ filtration has \emph{height} $n$. Experience shows that the complexity of the height $n$ part grows exponentially as $n$ increases linearly. For example for the sphere $S_{(p)}$, the height 0 elements are exactly the elements in $\pi_0$, the height 1 elements correspond to $\im j$, and the height 2 elements are associated with \emph{tmf}. Following \cite{RavenelOrange}, determining the chromatic filtration of an element $f \in \pi_{\ast}(X)$ can be done by a recursive algorithm that is based on the following steps.
\begin{enumerate}[leftmargin=*,labelindent=2em]
\item[Step 1:] Find a non-nilpotent self-map $v$ on $X$, not necessarily of degree 0. 
\item[Step 2:] Since $v$ acts on $\pi_{\ast}(X)$ by post-composition, we can use it to break $\pi_{\ast}(X)$ into a $v$-periodic part and a (power of) $v$-torsion part.
\item[Step 3:] If the element $f$ is $v$-periodic, then it will be detected in some cohomology theory, and we are done. If not, then it lifts to the cofiber of some power of $v$, and we repeat this process by replacing $X$ with this cofiber. This process increases the height of $f$ by 1.
\end{enumerate}
The execution and good behavior of this algorithm require the following ingredients:
\begin{enumerate}
\item the existence of a non-nilpotent self-map $v$ on every finite complex,
\item some sort of uniqueness for such a self-map $v$,
\item computable cohomology theories detecting $v$-periodicity.
\end{enumerate}
The first two points are exactly the content of one of the deepest theorems in chromatic homotopy theory called the \emph{Periodicity Theorem} of Devinatz-Hopkins-Smith. Given any finite complex $X$, this theorem says that there is an essentially unique non-nilpotent self-map on $X$ called a \emph{$v_n$-self-map}. Finally, the last point is taken care of by the existence of the \emph{Morava $K$-theories} $K(n)$ \cite{MoravaKth} which are field spectra\footnote{and thus admit a K\"unneth isomorphism, i.e., are computable.} and detect exactly $v_n$-periodicity.

Even over very nice base schemes (for example algebraically closed fields of characteristic 0), this algorithm does not apply in the category of motivic spectra. Motivic Morava $K$-theories detecting $v_n$-periodicity have been constructed in \cite{Borghesi} and \cite{HKOrem}, and even though they are not quite field spectra, they are computable over nice base schemes. The main issue comes from the lack of a periodicity theorem, and the fact that there is more periodicity to consider than just $v_n$-periodicity. In this paper we will explain this phenomenon in the easiest case: for 2-local cellular motivic spectra over $\Spec \C$. Denote by $S$ the 2-local motivic sphere spectrum over $\Spec \C$. The first step in the chromatic approach to compute $\piss(S)$ already fails, as there are two very different non-nilpotent self-maps
\begin{equation*}
S \stackrel{2}{\lto} S \qquad \text{ and } \qquad S^{1,1} \stackrel{\eta}{\lto} S.
\end{equation*}
This means that the process has to be refined, and that the linear ordering of periodicities has to be replaced by a more complex lattice. It turns out that in the same way that there is a $v_n$-periodicity story starting with $v_0=2$, Haynes Miller suggested that there could be a similar story starting with $w_0 = \eta$. The first evidence is provided by Andrews in \cite{Andrews}, where he also pins down the notation of \emph{$w_n$-periodicity}. Even though no precise definition is given, he shows that in the same way that $S/v_0$ admits a $v_1^4$-self-map, the motivic 2-cell complex $S/w_0$ also admits a $w_1^4$ self-map. The main goal of our paper is to give a precise definition of $w_n$-periodicity, as well as construct motivic field spectra $K(w_n)$ that detect $w_n$-periodicity. We refer to Theorem \ref{thm:kwnEoo} for more details.

The motivic fields $K(w_n)$ are also directly related to the structure of the category of finite 2-local motivic spectra. From the Periodicity Theorem, one can rewind back to one of the deepest and pioneering theorems in chromatic homotopy theory : the \emph{Thick Subcategory Theorem} proved by Devinatz-Hopkins-Smith in \cite{DHS}. This theorem is equivalent to the Nilpotence Theorem, and in fact implies the Periodicity Theorem. In some sense it is a more global way of understanding the chromatic filtration, without zooming in on a specific object. This theorem describes a filtration  
\begin{equation*}
\ast = \catC_{-1} \subsetneq \catC_0 \subsetneq \catC_1 \subsetneq \cdots \subsetneq \catC_{\infty} = \FinpSpt
\end{equation*}
by height, on the whole category $\FinpSpt$ of finite $p$-local spectra. The category $\catC_n$ is the subcategory of acyclics for the Morava $K$-theory spectrum $K(n)$ of height $n$. This filtration is exhaustive, Hausdorff, and admits no refinement by any thick subcategory. In fact, the subcategories $\catC_n$ turn out to further be prime ideals, i.e., this filtration gives a complete description of the \emph{Balmer spectrum of} $\FinpSpt$. Although it is not true in a general tensor triangulated category that all prime ideals come from field spectra, this happens to be the case for $\FinpSpt$. 

Motivically, the study of thick subcategories started with the work in \cite{Joa}. In that paper some thick subcategories of $\FinSpt_{\C}$ were constructed, with the feeling that this is a very hard problem. In fact, even the Balmer spectrum (i.e., just the thick prime ideals) of the category of finite $p$-local motivic spectra has not been computed over any base scheme. A slightly easier problem that we will consider is the Balmer spectrum of the category of \emph{cellular} finite $p$-local spectra. Denote this category (over the base $\Spec \C$) by $\FinCell_{\C}$. Since the spectra $K(w_n)$ are new motivic fields, it is now clear that the Balmer spectrum of $\FinCell_{\C}$ is more complicated than the one of $\FinSpt$. Morel showed \cite{Morelpi0} that if the base scheme is a perfect field $k$ of characteristic different than $2$, then
\begin{equation*}
\pi_{-(\ast,\ast)} \Ss \cong \KMW(k)
\end{equation*}
is the Milnor-Witt $K$-theory of $k$, which contains the Grothendieck-Witt group $K^{\text{MW}}_0(k) = \GW(k)$ in degree 0. In \cite{BalmerSSS}, Balmer considers a natural map
\begin{equation*}
\rho \colon \Spc(\FinCell_{k}) \lto \Spec(\GW(k))
\end{equation*}
from the spectrum of finite motivic cellular spectra. This maps send a thick prime ideal $\mathfrak{p}$ to the prime ideal of elements $f \in \GW(k) = [\Ss,\Ss]$ such that $Cf \notin \mathfrak{p}$, where $Cf$ denotes the 2-cell complex given by the cofiber of the map $f$. He shows by general methods that in this case, the map $\rho$ is surjective. In the case when $k$ is either $\C$ or a finite field $\F_q$, the only non-trivial prime ideals of $\GW(k)$ are given by $(p)$, and the surjectivity of $\rho$ was already known since the motivic $v_n$-story covers these ideals. In \cite{Kelly} for the case of finite fields, and more generally in \cite{HO}, it is shown that Balmer's map $\rho$ factors further through the surjective map $\rho^{\bullet}$ as in
\begin{center}
\begin{tikzpicture}
\matrix (m) [matrix of math nodes, row sep=3em, column sep=4em]
{ \Spc(\FinCell_{k}) & \Spech(\KMW(k)) \\
                    & \Spec(\GW(k)), \\};
\path[thick, -stealth]
(m-1-1) edge node[above] {$ \rho^{\bullet} $} (m-1-2)
(m-1-1) edge node[below = 4pt, left=4pt] {$ \rho $} (m-2-2)
(m-1-2) edge  (m-2-2);
\end{tikzpicture}
\end{center}
where $\Spec^{\text{h}}$ denotes the space of homogeneous prime ideals. The space $\Spech(\KMW(k))$ has been computed in \cite{Thornton}, the remaining task for understanding the thick prime ideals of finite cellular motivic spectra is thus to identify the fibers of this map. However, this is no easy task, even in what is considered to be the easiest case, i.e., over $\Spec \C$. In this case, the lattice of homogeneous prime ideals of $\Spech(\KMW(\C))$ is given by 
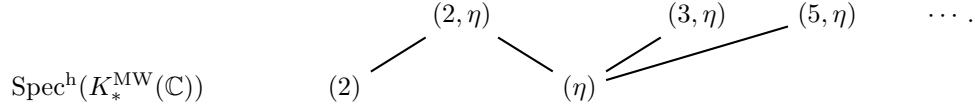
\begin{center}
\begin{tikzpicture}
\matrix (m) [matrix of math nodes, row sep=1em, column sep=2em]
{ & &    & (2,\eta) &        & (3, \eta) & (5,\eta) & \cdots. \\
  \Spech(\KMW(\C)) & & (2) &          & (\eta) & & &  \\};
\path[thick]
(m-2-3) edge (m-1-4)
(m-2-5) edge (m-1-4)
(m-2-5) edge (m-1-6)
(m-2-5) edge (m-1-7);
\end{tikzpicture}
\end{center}
Some explicit thick prime ideals have been constructed in \cite{HO} in the case when the base field $k$ admits an embedding $k \inj \C$. For simplicity, let's from now on only work over $\Spec \C$. Recall that there is an adjunction
\begin{equation*}
\Real \colon \SptC \adj \Spt \colon \Sing,
\end{equation*}
where $\Real$ is the \emph{Betti realization} functor induced by taking $\C$-points. In \cite{HO}, Heller-Ormsby show that the thick prime ideal $\Real^{-1}(\catC_n)$ is the subcategory of acyclics for the motivic field spectrum $\Sing(K(n))$. The thick prime ideal generated by $\Sing(K(0))$ sits over $(\eta)$, while $\Sing(K(n))$ for $n > 0$ sits over $(p,\eta)$ for the appropriate prime, very much like the classical picture except that $(0)$ becomes $(\eta)$. In the cellular case, being a motivic field is a weaker condition, since in particular it is implied by the coefficients being a graded field. In this case, the paper \cite{HO} constructs a motivic field from the motivic spectrum $KT$ of \cite{HornbostelHermitian} representing \emph{(higher) Witt groups}. This cellular motivic field generates another thick prime ideal and lives over the ideal $(2)$.

In Theorem \ref{thm:kwnEoo} we construct more cellular motivic fields $K(w_n)$ for every $n \in \mathbb{N}_0$. Because $w_0 = \eta$, the motivic spectrum $K(w_0)$ agrees on homotopy groups with the cellular field of \cite{HO}. These new motivic fields sit above $\Spech(\KMW(\C))$ as is shown in the diagram
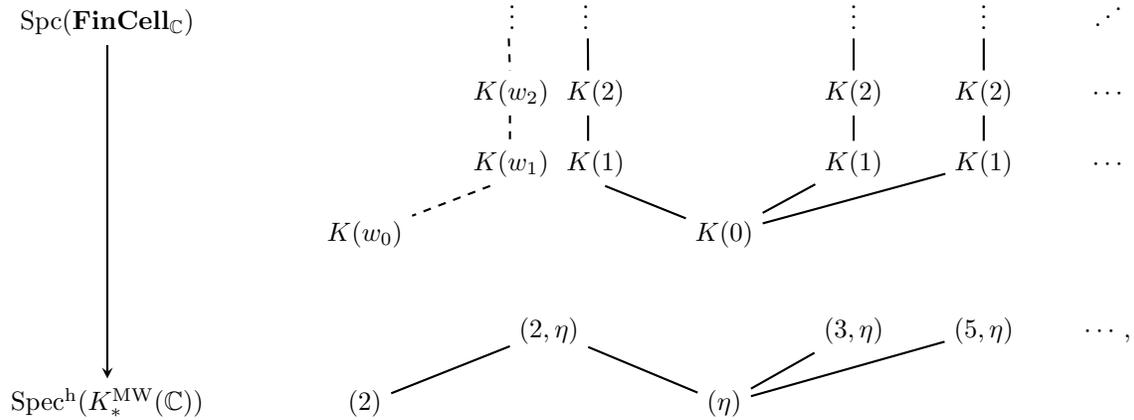
\begin{center}
\begin{tikzpicture}
\matrix (m) [matrix of math nodes, row sep=1em, column sep=2em]
{ \Spc( \FinCell_{\C}) & &    & \vdots \ \ \ \ \ \ \ \vdots &        & \vdots & \vdots & \iddots \\
                  & &    & K(w_2) \ \  K(2) &        & K(2) & K(2) & \cdots \\                
                  & &    & K(w_1) \ \  K(1) &        & K(1) & K(1) & \cdots \\                
                  & & K(w_0) &  & K(0)       &   &   &   \\                
                  & &    & &        & &  &  \\  
   				  & &    & (2,\eta) &        & (3, \eta) & (5,\eta) & \cdots, \\
  \Spech(\KMW(\C))& & (2)&          & (\eta) & & &  \\};
\path[thick]
(m-3-4) edge (m-4-5)
(m-2-4.330) edge (m-3-4.30)
(m-1-4.323.5) edge (m-2-4.30)

(m-3-6) edge (m-4-5)
(m-2-6) edge (m-3-6)
(m-1-6) edge (m-2-6)

(m-3-7) edge (m-4-5)
(m-2-7) edge (m-3-7)
(m-1-7) edge (m-2-7)

(m-7-3) edge (m-6-4)
(m-7-5) edge (m-6-4)
(m-7-5) edge (m-6-6)
(m-7-5) edge (m-6-7);
\path[thick, -stealth]
(m-1-1) edge (m-7-1);
\path[thick, dashed]
(m-3-4) edge (m-4-3)
(m-2-4.210) edge (m-3-4.150)
(m-1-4.216) edge (m-2-4.150);
\end{tikzpicture}
\end{center}
where the symbol $K(n)$ stands for the motivic field $\Sing(K(n))$. From this picture, it is tempting to conjecture that the containment of the thick ideals generated by $K(w_n)$ is similar to the $K(n)$ story (the above dotted conjectural lines represent this containment), and that these fields only exist at $p=2$. In any case, our methods employed to detect $w_n$-periodicity and to construct $K(w_n)$ in Section \ref{sec:kwn} do not generalize in an obvious way to odd primes. Work in progress of Barthel-Heard-Krause investigate such fields at odd primes.

Finally, after constructing motivic fields $K(w_n)$ that detect $w_n$-periodicity and their connective versions $k(w_n)$, it is desirable to also have a more global motivic spectrum containing all the $w_i$. In classical chromatic homotopy theory, the Brown-Peterson spectrum $BP$ is necessary to enunciate the Nilpotence Theorem. This theorem states that given any $p$-local finite complex $X$, a self-map $\Sigma^{\ast} X \lto X$ is nilpotent if and only if it is nilpotent in $BP$-homology. For example in the case of a single cell $X = S_{(p)}$, this theorem recovers Nishida's Theorem since $BP_{\ast}$ is torsion-free. Motivically, it is easy to see that the natural motivic analogue $BPGL$ does not detect nilpotence. We already pointed out that the Hopf map $\eta \colon S^{1,1} \lto S$ is not nilpotent. Moreover, by construction, the algebraic cobordism spectrum $MGL$ does not detect $\eta$, and thus neither does the spectrum $BPGL$. In Section \ref{sec:wBP} we construct an $E_{\infty}$ motivic ring spectrum $\wBP$ with homotopy groups
\begin{equation*}
\piss(\wBP) \cong \F_2[w_0, w_1, \ldots],
\end{equation*}
with the hope that $BPGL$ together with $wBP$ detect nilpotence. Unfortunately, this turns out to not be the case. More precisely, there is an element $d_1 \in \pi_{32,18}({\widehat{S}}_2)$ which is non-nilpotent by work of \cite{StableStems}, and which is not detected by either $BPGL$ or $\wBP$. At this point we should mention current work in progress of Barthel-Heard-Krause, which organizes $w_n$-periodicity in a bigger framework. Their idea is to rewind back to the $E_1$-page of the May spectral sequence, and consider all the periodicity that can occur from May's elements $h_{ij}$. In their work, they constructed a cellular motivic ring spectrum that detects the element $d_1$, as well as more motivic spectra (which are neither ring nor fields) detecting more $h_{ij}$-periodicity. One can hope that their additional spectra provide all the tools to detect nilpotence in the motivic setting.

\subsection{Organization}\mbox{}

\vspace{0.2cm}

\noindent \textbf{Section \ref{sec:Prereq}.} In this section we first describe the setting in which we work, which is the category of cellular $\Ct$-modules. This includes the following : computing the relevant Steenrod algebra and its dual, deriving some important properties, and setting up an appropriate Adams spectral sequence.

\vspace{0.3cm}

\noindent \textbf{Section \ref{sec:kwn}.} This section contains the construction of the motivic fields $K(w_n)$. This goes through first constructing connective versions $\kwn$, endowing them with an $E_{\infty}$ ring structure, and finally inverting multiplication by $w_n$.

\vspace{0.3cm}

\noindent \textbf{Section \ref{sec:wBP}.} This section contains the construction of the motivic spectrum $\wBP$, and its truncations $\wBPn$.

\subsection{Acknowledgment}
 
The author is grateful for contributions from Dan Isaksen, Nicolas Ricka, Jens Hornbostel, J.D. Quigley, Mark Behrens, Prasit Bhattacharya, Dominic Culver and Tyler Lawson.

\section{The Category of $\Ct$-modules and its Steenrod algebra} \label{sec:Prereq}

In this section we will describe the general framework in which all spectra will be constructed. We first start by recalling notation in \ref{subsec:motspectra} regarding the motivic setting, as well as set up the category $\CtCell$ of cellular (2-completed) $\Ct$-modules. In \ref{subsec:CtHHandSteenrod} we describe the $\Ct$-induced mod 2 Eilenberg-Maclane spectrum $\HH$, and its Steenrod algebra of operations and co-operations. This spectrum plays in $\CtCell$ the role that $\HF$ plays in $\SptC$, and will serve as a building block for the Postnikov tower constructions of Section \ref{sec:kwn}. In \ref{subsec:Pst} we further study this Steenrod algebra and give the relevant definitions which lead to the definition of $w_n$ periodicity. Finally, in \ref{subsec:HHAdamsss} we briefly describe the $\HH$-based Adams spectral sequence in the category $\CtCell$. This spectral sequence will be used in several places in Sections \ref{sec:kwn} and \ref{sec:wBP}.

\subsection{Cellular motivic spectra and $\Ct$-modules} \label{subsec:motspectra}

Denote by $\SptC$ the category of \emph{motivic spectra over $\Spec \C$} as introduced in \cite{MorelVoevodsky}. It is a closed symmetric monoidal category with smash product $- \smas -$, and internal motivic function spectrum denoted by $F(-,-)$. We will denote the 1-dimensional motivic spheres by
\begin{itemize}
\item $S^{1,0}$ for the \emph{simplicial sphere}, i.e., the constant presheaf on the simplicial circle $\Delta^1/\partial \Delta^1$,
\item $S^{1,1}$ for the \emph{geometric sphere}, i.e., the suspension spectrum of the group scheme\footnote{which we could call the punctured affine line, since we will not use its group structure.} $\G_m$.
\end{itemize}
Smashing together various copies of these spheres gives a bigraded family of spheres $S^{s,w}$ for $s,w \in \Z$, where $s$ is called the \emph{stem} and $w$ is called the \emph{weight}. Denote the bigraded homotopy groups of a motivic spectrum $X$ by $\piss(X)$, where the homotopy group $\pisw(X)$ in stem $s$ and weight $w$ is given by homotopy classes of maps from $S^{s,w}$ into $X$. Recall from \cite{DuggerIsaksen} that a motivic spectrum is called \emph{cellular} if it can be built out of spheres $S^{s,w}$ under filtered colimits. Denote by $\SCell$ the category of \emph{cellular motivic spectra over $\Spec \C$}, constructed as the right Bousfield localization at the set of spheres $\{ S^{s,w} \}_{s,w \in \Z}$. The weak equivalences in $\SCell$ are thus given by $\piss$-isomorphisms. The Bousfield localization is part of an adjunction
\begin{equation*}
\SCell \overunder{}{C}{\adjinc} \SptC,
\end{equation*}
where the unit is a weak equivalence in $\SCell$, and the counit $C(X) \lto X$ is a $\piss$-isomorphism. This discussion can be carried out both in the world of 
\begin{itemize}
\item presentable, closed symmetric monoidal $\infty$-categories, following \cite{Robalo},
\item cellular, closed symmetric monoidal model categories, via the motivic symmetric spectra of \cite{Jardine}, \cite{Pelaez}, and the theory of right Bousfield localization in that setting following \cite{BarnesRoit}.
\end{itemize}
In this paper we will be working in $\SCell$ as all spectra constructed will be cellular. This has in particular the advantages that a motivic spectrum $X \in \SCell$ is contractible if and only if its homotopy groups $\piss(X)$ vanish, and that our spectral sequences converge.

In \cite{VoeZ2}, Veovodsky introduced the motivic mod 2 Eilenberg-Maclane spectrum $\HF \in \SptC$, and computed its coefficients 
\begin{equation*}
\piss(\HF) \cong \F_2[\tau].
\end{equation*}
The polynomial class $\tau$ is in bidegree $(0,-1)$, and is related to the Tate twist. In the further papers \cite{Voered}, \cite{VoeEM}, Voevodsky computed the Hopf algebra structure of the motivic mod 2 Steenrod algebra and its dual, with algebra structure on the dual given by
\begin{equation} \label{eq:dualSteenrod}
\ACdual \cong \quotient{\F_2[\tau][\xi_1, \xi_2, \ldots][\tau_0, \tau_1, \ldots]}{\tau_i^2 = \tau \xi_{i+1}}.
\end{equation}
The element $\tau$ is primitive, and the coproduct on the $\xi_i$'s and $\tau_i$'s is given by the same formula as the classical odd primary formula. In \cite{HKOrem}, \cite{HKOAdams} it is shown that the motivic spectrum $\HF$ is in fact cellular. One can thus consider the $\HF$-based motivic Adams spectral sequence as constructed in \cite{MorelAdams}, \cite{DuggerIsaksenMASS}, \cite{HKOAdams}, which computes the homotopy groups of the 2-completed sphere $(\Ss)^{\smas}_2$. For degree reasons, multiplication by $\tau$ on $\HFco(\Ss)$ survives to the $E_{\infty}$-page, producing a map
\begin{equation*}
S^{0,-1} \stackrel{\tau}{\lto} (\Ss)^{\smas}_2.
\end{equation*}
This map exists after completion at any prime $p$, but does not lift to an integral map $S^{0,-1} \lto \Ss$. 

\begin{remark}[2-completion]
Since $\tau$ does not exist before completing the target, we will from now on exclusively work in the 2-completed category without further mention, and drop the completion symbol from the notation. For example, we will denote the category of 2-completed cellular motivic spectra by $\SCell$, the 2-completed motivic sphere spectrum by $\Ss$, the smash product over the 2-completed sphere by $- \smas -$, etc. With this
notation, the motivic Adams spectral sequence produces a non-trivial map $S^{0,-1} \stackrel{\tau}{\lto} \Ss$ that can be seen as an element in the homotopy group $\pi_{0,-1}\Ss$. 
\end{remark}
Denote the cofiber of $\tau \in \pi_{0,-1} \Ss$ by
\begin{equation*}
S^{0,-1} \stackrel{\tau}{\lto} \Ss \stackrel{i}{\lto} \Ct.
\end{equation*}
It is shown in \cite{GheCt} that the motivic 2-cell complex $\Ct$ admits a unique $E_{\infty}$ ring structure, producing thus a closed symmetric monoidal category $(\CtMod, - \smas_{\Ct} -)$ of $\Ct$-modules. The usual adjunction from the ring map $\Ss \lto \Ct$ restricts to an adjunction 
\begin{equation*}
\SCell \stackrel{- \smas \Ct}{\adj} \CtCell
\end{equation*}
on cellular objects. From now on, we will call an object $X \in \CtCell$ a $\Ct$-module, omitting the word cellular.

\begin{remark}[Working in $\Ct$-modules] \label{rem:workinCt}
Similarly with the fact that the Morava $K$-theories $K(n)$ are 2-completed\footnote{except for $K(0) \simeq H\Q$.}, i.e., are $(\Ss)^{\smas}_2$-modules, the spectra $\Kwn$ are naturally $\Ct$-modules. In the case of $K(w_0)$, this can be seen from the relation $0 = \tau \eta^4 \in \pi_{4,3} (\widehat{S}_2)$. Since $K(w_0)$ contains $\eta^{-1}$, this forces $\tau$ to act by zero on it, which in this case is sufficient to promote a motivic spectrum to a $\Ct$-module. We will therefore work in the category of $\CtCell$ in which we will construct the motivic fields $K(w_n)$.
\end{remark}

\subsection{$\Ct$-linear $\HH$-(co)homology and its (co)operations} \label{subsec:CtHHandSteenrod}

We will now exclusively be working in $\CtCell$, i.e., with cellular $\Ct$-modules and with $\Ct$-linear maps between them. Most invariants of the underlying motivic spectrum of a $\Ct$-module $X$ can be rewritten in this category. For example, the usual adjunction describes its homotopy groups by
\begin{equation*}
\pisw(X) \cong \left[ \Sigma^{s,w} \Ct, X \right]_{\Ct}.
\end{equation*}
The analog of the mod 2 Eilenberg-Maclane spectrum in this category is the $\Ct$-induced Eilenberg-Maclane spectrum $\HH \coloneqq \HF \smas \Ct$. As we explain in Remark \ref{rem:CtlinHHhom}, the following definition recovers the mod 2 (co)homology of the underlying motivic spectrum.

\begin{defn}
Given a $\Ct$-module $X$, define its \emph{$\Ct$-linear $\HH$-homology} to be
\begin{equation*}
\HHho(X) \coloneqq \piss( \HH \smas_{\Ct} X ).
\end{equation*}
This is naturally coacted upon by the \emph{$\Ct$-linear co-operations of $\HH$} defined by $\Aaadual \coloneqq \piss(\HH \smas_{\Ct} \HH)$. Similarly define the \emph{$\Ct$-linear $\HH$-cohomology} of $X$ to be
\begin{equation*}
\HHco(X) \coloneqq \pi_{-\ast,-\ast}( F_{\Ct}(X , \Sigma^{-1,1} \HH) ),
\end{equation*}
where the reason for the shift is apparent in Remark \ref{rem:CtlinHHhom}. This is naturally acted upon by the \emph{$\Ct$-linear operations of $\HH$} defined by $\Aaa \coloneqq \pi_{-\ast,-\ast}F_{\Ct}(\HH,\HH)$. We will sometimes call $\Aaa$ (respectively $\Aaadual$) the ($\Ct$-linear) (respectively \emph{dual})\emph{ Steenrod algebra of $\HH$}.
\end{defn}

\begin{remark}[$\Ct$-linear $\HH$-(co)homology versus $\HF$-(co)homology] \label{rem:CtlinHHhom}
The $\Ct$-linear $\HH$-homology of $X$ is naturally isomorphic to the $\HF$-homology of its underlying motivic spectrum since
\begin{equation*}
\HHho(X) = \piss(\HH \smas_{\Ct} X) \cong \piss(\HF \smas \Ct \smas_{\Ct} X) \cong \piss( \HF \smas X) =\HFho(X).
\end{equation*}
The inclusion of the bottom cell $\Ss \lto \Ct$ induces a map $\ACdual \lto \Aaadual$ from the $\HF$ dual Steenrod algebra to the $\Ct$-linear dual $\HH$-Steenrod algebra. This means that given a $\Ct$-module $X$, the coaction of $\ACdual$ factors through the coaction of $\Aaadual$. Similarly in cohomology, the forget-hom adjunction gives the equivalence
\begin{equation*}
\HHco(X) = \pi_{-\ast,-\ast}( F_{\Ct}(X , \Sigma^{-1,1} \HH) ) \cong \pi_{-\ast,-\ast}( F_{\Ct}(X, F(\Ct,\HF)) ) \cong \pi_{-\ast,-\ast}( F(X,\HF) ) \cong \HFco(X).
\end{equation*}
The inclusion of the bottom cell $\Ss \lto \Ct$ again produces a map $\AC \lto \Aaa$, and so the action of $\AC$ factors through the action of $\Aaa$. 
\end{remark}

The computations of $\Aaa$ and $\Aaadual$ follow easily from Voevodsky's computation of $\AC$ and $\ACdual$.

\begin{prop} [{\cite[Remark 5.6]{GheCt}}] \label{prop:HHSteenrod}
The $\Ct$-linear co-operations of $\HH$ are given by the Hopf algebra
\begin{equation*}
\Aaadual \cong \F_2[\xi_1, \xi_2, \ldots] \otimes E(\tau_0, \tau_1, \ldots),
\end{equation*}
with bidegrees given by $|\xi_n| = (2^{n+1}-2, 2^n -1)$ and $|\tau_n| = (2^{n+1}-1, 2^n -1)$, and coproduct
\begin{align*}
\Delta(\xi_n) &= \sum_{i = 0}^n \xi_{n-i}^{2^i} \otimes \xi_i =  \xi_n \otimes 1 + \xi_{n-1}^2 \otimes \xi_1 + \cdots + \xi_{n-i}^{2^i} \otimes \xi_i + \cdots + 1 \otimes \xi_n, \\ 
\Delta(\tau_n) &= \tau_n \otimes 1 + \sum_{i = 0}^n \xi_{n-i}^{2^i} \otimes \tau_i =  \tau_n \otimes 1 + \xi_n \otimes \tau_0 + \xi_{n-1}^2 \otimes \tau_1 + \cdots + \xi_{n-i}^{2^i} \otimes \tau_i + \cdots + 1 \otimes \tau_n. \\ 
\end{align*}
\end{prop}

The advantage of working with the coaction of $\Aaadual$ instead of $\ACdual$ is now apparent by comparing Proposition \ref{prop:HHSteenrod} with equation \eqref{eq:dualSteenrod}. First, $\Aaadual$ is smaller and more regular, which will be convenient for computations. Second, since $\tau \in \pi_{0,-1} \Ss$ is nullhomotopic on any $\Ct$-module $X$, it will act as zero on any algebraic invariant of $X$. Morally, it is thus natural to expect that the coaction of $\ACdual$ on the homology $\HFho(X)$ should factor through the quotient $\quotient{\ACdual}{\tau} \cong \Aaadual$. Working with $\Ct$-linear $\HH$-homology is a way of making this remark precise. The exact same remark applies to cohomology, where the computation of the $\Ct$-linear Steenrod algebra $\Aaa$ follows from \cite[Proposition 5.5]{GheCt} and is given by the quotient $\Aaa \cong \quotient{\AC}{\tau}$.

\begin{conv}
Given a $\Ct$-module $X$, we will always consider its $\Ct$-linear $\HH$-homology (i.e., its $\HF$-homology) endowed with the coaction of $\Aaadual$. Similarly, its $\Ct$-linear $\HH$-cohomology will always be considered as an $\Aaa$-module.
\end{conv}

To state another crucial advantage of $\Aaa$ over $\AC$ we need the following definition.

\begin{defn}[Chow degree]
Let $A_{\ast,\ast}$ be a bigraded abelian group. The \emph{Chow degree} of an element $x \in A_{s,w}$ is given by the difference $s - 2w$. The bigraded group $A_{\ast,\ast}$ splits as a sum of its summands in a fixed Chow degree, which ranges through $\Z$.
\end{defn}

\begin{remark}[The Chow degree on motivic Steenrod algebras] \label{rem:ChowSteerodalg}
Both $\AC$ and $\Aaa$ are generated as algebras (over $\F_2[\tau]$ and $\F_2$ respectively) by the Steenrod squares $\Sqnone$. The even squares $\Sqn$ are in Chow degree 0, while the odd squares $\Sqnn$ are in Chow degree 1. It follows that the whole Steenrod algebra $\Aaa$ is concentrated in non-negative Chow degrees, i.e., is bounded below by 0. This in particular allows recursive arguments on the Chow degree. On the other side, since in this cohomological setting $|\tau| = (0,1)$ is in Chow degree $-2$, the Steenrod algebra $\AC$ of Voevodsky is non-vanishing in all Chow degrees and does not allow recursive arguments of this type.
\end{remark}

Finally, let's mention that $\HHco$-cohomology satisfies a K\"unneth formula.

\begin{prop} \label{prop:Kunneth}
Given two $\Ct$-modules $X$ and $Y$, there is a K\"unneth isomorphism
\begin{equation*}
\HHco(X \smas Y) \cong \HHco(X) \otimes_{\F_2} \HHco(Y)
\end{equation*}
of $\Aaa$-modules, where the right hand side has the diagonal $\Aaa$-module structure. 
\begin{proof}
Consider the motivic K\"unneth spectral sequence from \cite{DuggerIsaksen}, \cite{Totaro}
\begin{equation*}
\Tor^{s,t,w}_{\HHco} \big( \HHco(X), \HHco(Y) \big) \Longrightarrow \HH^{t-s,w}(X \smas Y),
\end{equation*}
where $s$ is the homological degree, and $(t,w)$ are the usual two internal degrees. Since $\HHco \cong \F_2$ is a field, this spectral sequence is concentrated in homological degree $s=0$ and thus collapses, giving the desired result.
\end{proof}
\end{prop}

\subsection{The Motivic Margolis elements $\Pst$} \label{subsec:Pst}

In this section we will set-up some notation and formulas in the Steenrod algebra $\Aaa$. These formulas are the motivic adaptation of classical formulas in $\Acl$ proven by Milnor, see for example \cite[Section 2]{BP}. There are two different ways of showing these formulas in the motivic setting
\begin{enumerate}
\item either by brute-force, by adapting the classical proof to the motivic setting, or
\item by transporting them via a map between the classical and motivic Steenrod algebras.
\end{enumerate}
We chose to use the second option. Consider the injective map
\begin{equation} \label{eq:mapofSteenrod}
\Acl \inj \Aaa
\end{equation}
of Hopf algebras, defined in \cite[Section 2.1.3]{StableStems}, where $\Acl$ denotes the mod 2 classical Steenrod algebra. One can for example define it as the dual to the natural quotient map between dual Steenrod algebras
\begin{equation} \label{eq:mapofdualSteenrod}
\Aaadual = \F_2[\xi_1, \xi_2, \ldots] \otimes E(\tau_0, \tau_1, \ldots) \lto \Acldual = \F_2[\xi_1, \xi_2, \ldots].
\end{equation}
Observe that a motivic Adams-Novikov version of this map is used in the detection part of \cite[Definition 2.5]{Andrews}.

\begin{remark} \label{rem:gradedbyweight}
It is easy to see that both maps are graded if the motivic bigraded object is consider as simply graded by the weight. Restricting to the weight in the motivic setting feels artificial, but turns out to be useful for the following reason. Denote by $c(-)$ the conjugation map on both the classical and motivic Steenrod algebras (and their duals). By analyzing the coproduct on $\Aaadual$, it is easy to see that the ring map \eqref{eq:mapofdualSteenrod} is in fact a \emph{graded} map of bialgebras over $\F_2$. Since both are connected Hopf algebras, the conjugation $c(-)$ is uniquely determined, and thus both maps \eqref{eq:mapofSteenrod} and \eqref{eq:mapofdualSteenrod} are maps of Hopf algebras.
\end{remark}

\begin{nota}[Margolis' $\Pst$] \label{nota:Pst}
Denote by $\Pst \in \Aaa$ the element dual to $\xi_t^{2^s} \in \Aaadual$, by dualizing in the canonical monomial basis. In the case $s=0$, we will simply denote this element by $P_t = P^0_t$.
\end{nota}

\begin{example}
Since the motivic Steenrod algebra at $p=2$ admits a slightly different notation than the classical one, let's look at low dimensional elements. When $s=0$, these elements are the sequence $P_1 = \Sqq$, $P_2 = [\Sqq, \Sqqqq]$, etc, and they follow the same pattern as Milnor's $Q_t$ sequence. Observe that $Q_0 = \Sq$ does not appear in this notation, and the sequence starts with $P_1 = \Sqq$. The sequence $\{ P_t \}$ can thus be seen as a doubled version of the classical Milnor sequence $\{ Q_t \}$.
\end{example}

\begin{lemma} \label{lem:Pstext}
The element $\Pst \in \Aaa$ is exterior if and only if $s < t$. Moreover, the subalgebra generated by the elements $\Pt$ is an exterior commutative algebra.
\begin{proof}
Let's also denote by $\Pst \in \Acl$ the classical element dual to $\xi_t^{2^s}$. Since the map \eqref{eq:mapofdualSteenrod} sends $\xi_t^{2^s}$ to $\xi_t^{2^s}$, its dual map \eqref{eq:mapofSteenrod} sends $\Pst$ to $\Pst$. It is proven in \cite[Lemma 15.1.4]{Margolis} that the classical $\Pst$'s are exterior if and only if $s < t$. The if part follows immediately and the only if part follows by injectivity of the map \eqref{eq:mapofSteenrod}.

In the classical setting, recall that the dual element to $\xi_t \in \Acldual$ is the Milnor primitive $Q_{t-1} \in \Acl$. It follows that the map \eqref{eq:mapofSteenrod} sends $Q_{t-1}$ to $\Pt$. Since the $Q_t$'s commute, then so do the $\Pt$'s, finishing the proof.
\end{proof}
\end{lemma}

\begin{nota}
Denote by $E(\Pt)$ the exterior algebra (in $\Aaa$) generated by $\Pt$, and by $\EPinf$ the exterior algebra (in $\Aaa$) generated by $P_1, P_2, \ldots$.
\end{nota}

Since $\Pt$ is exterior, one can consider \emph{Margolis homology with respect to $\Pt$}. Recall that given an $\Aaa$-module $M$, this is defined as the homology of the complex
\begin{equation*}
M \stackrel{\cdot \Pt}{\ltoo} M \stackrel{\cdot \Pt}{\ltoo} M,
\end{equation*}
i.e., by the formula $H(M;\Pt) = \ker \Pt / \im \Pt$. If $H(M;\Pt) = 0$, one says that $\Pt$ is \emph{exact }on $M$.

\begin{cor} \label{cor:Pstprimextexact}
For every $t$, the element $\Pt$ is primitive, exterior and exact on $\Aaa$.
\begin{proof}
Notice that $\Pt$ is primitive since it is dual to the indecomposable element $\xi_t$, and that it is exterior by Lemma \ref{lem:Pstext} applied with $s=0$. 

To show the vanishing of the $\Pt$-Margolis homology on $\Aaa$, we use the same strategy as in \cite[Proposition 19.1.1]{Margolis}. First of all, it is easy to see by inspection that the subalgebra $\EPinf$  has no Margolis homology for every $\Pt$. By a theorem of Milnor-Moore \cite{MilnorMoore}, the Hopf algebra $\Aaa$ is free over $\EPinf$, i.e., can be written as a direct sum $\oplus \EPinf$ as an $\EPinf$-module. It follows that $\Aaa$ has no $\Pt$-Margolis homology, i.e., that $\Pt$ is exact on $\Aaa$.
\end{proof}
\end{cor}

\begin{remark}[The Steenrod algebra and $w_n$-periodicity]
There is a tight relation between the periodic operator $v_n$ and the cohomology operation $Q_n$. This can for example be seen by the interplay between the homotopy and cohomology of the connective Morava $K$-theory spectrum $k(n)$. In the classical setting, these invariants are
\begin{equation*}
\pi_{\ast}(k(n)) \cong \F_2[v_n] \qquad \text{ and } \qquad H^{\ast}(k(n);\F_2) \cong \Acl // E(Q_n).
\end{equation*}
This can equivalently be seen in the Postnikov tower of $k(n)$, whose layers are given by Eilenberg-Maclane spectra $\HF$, which are attached via $Q_n$. The same relation exists motivically for the motivic Morava $K$-theories. An intuition for $w_n$-periodicity is that the relation between $w_n$ and Margolis' $\Pnn$ is the exact same as the relation between $v_n$ and $Q_n$.
\end{remark}

\subsection{The $\HH$-based $\Ct$-linear motivic Adams spectral sequence} \label{subsec:HHAdamsss}

In what follows, we will construct an $\HH$-based Adams spectral sequence in the category of $\Ct$-modules. For convergence issues, we will restrict to spectra $X$ that are bounded below in the sense that $\pi_{s,w} X = 0$ if $s < 0$. One can set it up as in \cite[Section 7]{DuggerIsaksenMASS} by replacing motivic spectra with $\Ct$-modules and $\HF$ with $\HH$. In this setting, an Adams resolution of a bounded below $\Ct$-module $X$ is given by a diagram of $\Ct$-modules
\begin{center}
\begin{tikzpicture}
\matrix (m) [matrix of math nodes, row sep=3em, column sep=5em]
{  X & X_1 & X_2 & \cdots \\
   K_0 & K_1 & K_2, &  \\};
\path[thick, -stealth]
(m-1-1) edge node[left] {$ g_0 $} (m-2-1)
(m-1-2) edge node[left] {$ g_1 $}  (m-2-2)
(m-1-3) edge node[left] {$ g_2 $}  (m-2-3)
(m-1-2) edge node[above] {$ f_0 $} (m-1-1)
(m-1-3) edge node[above] {$ f_1 $} (m-1-2)
(m-1-4) edge node[above] {$ f_2 $} (m-1-3);
\end{tikzpicture}
\end{center}
where $K_i$ is a wedge of suspensions of $\HH$ and $f_i$ is zero in $\HH$-(co)homology. Consider the fiber sequence $F \stackrel{f}{\lto} \Ct \stackrel{g}{\lto} \HH$ of $\Ct$-modules, where the map $\Ct \lto \HH$ comes from the unit $\Ss \lto \HF$. Then we can form a canonical Adams resolution as usual by inductively setting $K_i \simeq X \smas_{\Ct} F^{\smas_{\Ct} i} \smas_{\Ct} \HH$ and $X_i \simeq X \smas_{\Ct} F^{\smas_{\Ct} i}$, where $f_i$ and $g_i$ are induced from $f$ and $g$.

In \cite[Section 7]{DuggerIsaksenMASS} Dugger-Isaksen define a category $\langle \Ss \rangle_{\HF}$ for which the motivic Adams spectral sequence converges. This is the full subcategory of $\SptC$ containing cellular spectra, which is also closed by smashing with $\HF$. This last condition was necessary at that time since it was not known that $\HF$ was cellular. It was later proved in \cite{HKOAdams} that the mod 2 motivic Eilenberg-Maclane spectrum $\HF$ is cellular, and thus that $\langle \Ss \rangle_{\HF}$ is just the category of motivic cellular spectra. 

By copying \cite[Section 7]{DuggerIsaksenMASS}, we get a $\Ct$-linear $\HH$-based motivic Adams spectral sequence with $E_2$-term given by
\begin{equation*}
E_2 \cong \Ext_{\Aaa}(\HHco(X), \HHco(\Ct) ) \cong \Ext_{\Aaa}(\HHco(X), \F_2 ).
\end{equation*}
It remains to study what the $E_{\infty}$-page computes. Again by \cite[Section 7]{DuggerIsaksenMASS}, the $E_{\infty}$-page computes the homotopy groups of the homotopy limit of the semi-cosimplicial motivic spectrum
\begin{equation} \label{eq:towerCtmod}
\begin{tikzcd}
\HH \smas_{\Ct} X
\arrow[r, shift left]
\arrow[r, shift right]
& \HH \smas_{\Ct} \HH \smas_{\Ct} X 
\arrow[r]
\arrow[r, shift left=2]
\arrow[r, shift right=2]
& \cdots 
\end{tikzcd}
\end{equation}
where all cofaces are induced from $\Ct \lto \HH$. To compute this homotopy limit one can compare it with the $\HF$-tower of the underlying motivic spectrum of $X$
\begin{equation} \label{eq:towerSmod}
\begin{tikzcd}
\HF \smas X
\arrow[r, shift left]
\arrow[r, shift right]
& \HF \smas \HF \smas X
\arrow[r]
\arrow[r, shift left=2]
\arrow[r, shift right=2]
& \cdots 
\end{tikzcd}
\end{equation}
which we know totalizes to $X$ by \cite[Section 7]{DuggerIsaksenMASS}, since $X$ is already 2-complete. Since $\HH \smas_{\Ct} \Ct \simeq \HF \smas \Ss$, there are level-wise weak equivalences between these two towers\footnote{more precisely between the underlying tower of \eqref{eq:towerCtmod} and the tower \eqref{eq:towerSmod}.} which commute with the coface maps. Since forgetting the $\Ct$-module structure is a right adjoint, it follows that the underlying motivic spectrum of the totalization of the tower \eqref{eq:towerCtmod} is also $X$. There is thus a convergent $\Ct$-linear $\HH$-based motivic Adams spectral sequence
\begin{equation} \label{eq:HHAss}
\Ext^{s,t,w}_{\Aaa}(\HHco(X), \F_2 ) \Longrightarrow \pi_{t-s,w}(X).
\end{equation}

\begin{remark}
Even though both towers \eqref{eq:towerCtmod} and \eqref{eq:towerSmod} have the same underlying motivic spectrum, they do not live in the same category and thus do not produce the same spectral sequence. They converge to the same object, but the $E_2$-term of the Adams spectral sequence coming from \eqref{eq:towerCtmod} is smaller and more computable.
\end{remark}


\section{The Motivic Fields $\Kwn$} \label{sec:kwn}

The goal of this section is to construct motivic fields $\Kwn$, that detect $w_n$ periodicity. As explained in Remark \ref{rem:workinCt}, these spectra will be constructed in the category $\CtCell$. We will thus from now on exclusively work in the category $\CtCell$ and denote the smash product over $\Ct$ simply by $- \smas -$, homotopy classes of $\Ct$-linear maps by $[-,-]$, work with $\Ct$-linear $\HH$-homology and cohomology, etc. The motivic spectrum $\Kwn$ should be a ring spectrum with homotopy groups given by
\begin{equation*}
\piss(\Kwn) \cong \F_2[w_n^{\pm 1}].
\end{equation*}
The strategy is to first construct a connective version $\kwn$ with homotopy groups 
\begin{equation*}
\piss(\kwn) \cong \F_2[w_n],
\end{equation*}
endow it with a ring structure, and finally invert multiplication by $w_n$ to get $\Kwn$. In \ref{subsec:kwn} we will construct $\kwn$ and show that it has the correct homotopy, and appropriate cohomology. The construction is done along an inverse tower, which can be seen as a Postnikov tower of $\kwn$ (in the stem direction, for example). From this tower one can easily compute the homotopy and cohomology of $\kwn$ by the associated spectral sequences. In \ref{subsec:Eookwn} we will again use this tower to construct a ring map $\kwn \smas \kwn \lto \kwn$. We will then use Robinson's obstruction theory to rigidify it to an $E_{\infty}$ ring structure, which allows the definition of $\Kwn$.

\subsection{The construction of $\kwn$} \label{subsec:kwn}

Fix an $n \in \mathbb{N}_0$ until the end of the section. We will now construct $\kwn$ via its Postnikov tower, in the category of $\Ct$-modules. Recall that $\kwn$ should be a motivic ring spectrum whose homotopy groups are given by the polynomial ring
\begin{equation*}
\piss(\kwn) \cong \F_2[w_n],
\end{equation*}
where $w_n$ is an element detected by the cohomology operation $\Pnn$. This suggests that $\kwn$ could be constructed via a tower whose layers are copies of $\HH$, each of which is attached by $\Pnn$. We will call this tower the Postnikov tower of $\kwn$, as the layers are becoming more and more connected (in both the stem, and the weight). We proceed to explain this construction after the following notation convention.

\begin{nota}
Denote by $\rn \coloneqq |\Pnn| = |\xi_{n+1}| = (2^{n+2} -2, 2^{n+1} -1)$.
\end{nota}

\begin{cons}[The construction  of $\kwn$] \label{cons:constkwn}
Start the bottom of the tower with the fiber sequence
\begin{center}
\begin{tikzpicture}
\matrix (m) [matrix of math nodes, row sep=3em, column sep=5em]
{  \HH & \kww0 &  \\
        & \ast & \Sigma^{1,0} \HH, \\};
\path[thick, -stealth]
(m-1-1.346) edge node[above] {$ i_{-1} = \id $} (m-1-2)
(m-1-2) edge node[right] {$ p_{-1} $} (m-2-2)
(m-2-2) edge node[above] {$ k_{-1} = 0 $} (m-2-3.182);
\end{tikzpicture}
\end{center}
where $k_{-1}$ denotes the $-1^{\text{st}}$ $k$-invariant. The Postnikov truncation $\kww0$ has thus homotopy groups $\piss( \kww0) \cong \piss(\HH) \cong \F_2$. We now want to attach the second copy of $\HH$ via the Steenrod operation $\Pnn$. This requires a $k$-invariant $k_0$ that restricts to the operation $\Pnn$ on $\HH$ as shown in the diagram
\begin{center}
\begin{tikzpicture}
\matrix (m) [matrix of math nodes, row sep=3em, column sep=5em]
{   \HH & \kww0 & \Sigma^{\rn} \HH \\
       & \ast & \Sigma^{1,0} \HH, \\};
\path[thick, -stealth]
(m-1-1) edge node[above] {$ \id $} (m-1-2)
(m-1-2) edge node[right=8pt, below=-2pt] {$ $} (m-2-2)
(m-2-2) edge node[above] {$ $} (m-2-3);
\path[dashed, -stealth]
(m-1-1) edge[bend right = 15] node[below=5pt, left =20pt] {$ \Pnn  $} (m-1-3);
\path[dotted, -stealth]
(m-1-2) edge node[above] {$ \exists ? \ k_{0} $} (m-1-3);
\end{tikzpicture}
\end{center}
where we denote the bidegree of $\Pnn$ by $\rn = |\Pnn|$. There is obviously a unique such filler up to homotopy, which is $k_0 = \Pnn$. The next step is to take the fiber of $k_0$ to get the next stage in the tower, which we denote by $\kww1$ as shown in
\begin{center}
\begin{tikzpicture}
\matrix (m) [matrix of math nodes, row sep=3em, column sep=5em]
{  \Sigma^{\rn - (1,0)} \HH & \kww1 & \Sigma^{2\rn - (1,0)} \HH \\
   \HH & \kww0 & \Sigma^{\rn} \HH \\
        & \ast & \Sigma^{1,0} \HH. \\};
\path[thick, -stealth]
(m-1-1) edge node[above] {$ i_0 $} (m-1-2)
(m-2-1) edge node[above] {$ \id $} (m-2-2)
(m-1-2) edge node[right=8pt, below=-2pt] {$ p_0 $} (m-2-2)
(m-2-2) edge node[right=8pt, below=-2pt] {$ $} (m-3-2)
(m-3-2) edge node[above] {$  $} (m-3-3)
(m-2-2) edge node[above] {$ k_{0} $} (m-2-3);
\path[dotted, -stealth]
(m-1-2) edge node[above] {$  \exists ? \  k_{1} $} (m-1-3);
\path[dashed, -stealth]
(m-1-1) edge[bend right = 15] node[below=5pt, left =20pt] {$ \Pnn  $} (m-1-3)
(m-2-1) edge[bend right = 15] node[below=5pt, left =20pt] {$ \Pnn  $} (m-2-3);
\end{tikzpicture}
\end{center}
To continue the process, we need the existence of a $k$-invariant $k_1$ that restricts to $\Pnn$ on $\Sigma^{\rn - (1,0)} \HH$. Equivalently, we are trying to extend $\Pnn$ to $\kww1$ as in the diagram
\begin{center}
\begin{tikzpicture}
\matrix (m) [matrix of math nodes, row sep=3em, column sep=3em]
{  \Sigma^{-1,0} \kww0 & \Sigma^{\rn - (1,0)} \HH & \kww1 &  \kww0  \\
                      &  &\Sigma^{2\rn - (1,0)} \HH. & \\};
\path[thick, -stealth]
(m-1-1) edge node[above] {$ k_0 $} (m-1-2)
(m-1-2) edge node[above] {$ i_0 $} (m-1-3)
(m-1-3) edge node[above] {$ p_0 $} (m-1-4)
(m-1-2) edge node[below=8pt, left] {$ \Pnn  $} (m-2-3);
\path[dotted, -stealth]
(m-1-3) edge node[right] {$  \exists ? \  k_{1} $} (m-2-3);
\end{tikzpicture}
\end{center}
It follows that a $k$-invariant $k_1$ exists if and only if the composite
\begin{equation*}
\Sigma^{-1,0} \kww0 \stackrel{k_0}{\lto} \Sigma^{\rn-(1,0)} \HH \stackrel{\Pnn}{\lto}  \Sigma^{2\rn-(1,0)} \HH
\end{equation*}
is nullhomotopic, which is equivalent to having the relation $\Pnn k_0 = 0$ in the cohomology $\HH^{2\rn}(\kww0)$. If such a $k$-invariant $k_1$ exists, then the difference of two such extensions would factor through the map $p_0$, and so one can alter $k_1$ up to homotopy by the group of maps in $\left[ \kww0, \Sigma^{2\rn-(1,0)} \HH \right]$ modulo $k_0$-divisibles. In order to keep the outline of the construction clear, we postpone the proof of the relation $\Pnn k_0 = 0$, and the fact that the moduli of such extensions is trivial to Proposition \ref{prop:existkwn}. These two facts show that there exists a unique $k$-invariant $k_1$ up to homotopy, so the homotopy type $\kww1$ is uniquely defined and we can canonically continue to build the tower. We can now repeat the process by taking the fiber of $k_1$ as shown in the diagram
\begin{center}
\begin{tikzpicture}
\matrix (m) [matrix of math nodes, row sep=3em, column sep=5em]
{  \Sigma^{2\rn - (2,0)} \HH & \kww2 & \Sigma^{3\rn - (2,0)} \HH \\
	\Sigma^{\rn - (1,0)} \HH & \kww1 & \Sigma^{2\rn - (1,0)} \HH \\
   \HH & \kww0 & \Sigma^{\rn} \HH \\
        & \ast & \Sigma^{1,0} \HH. \\};
\path[thick, -stealth]
(m-1-1) edge node[above] {$ i_1 $} (m-1-2)
(m-2-1) edge node[above] {$ i_0 $} (m-2-2)
(m-3-1) edge node[above] {$ \id $} (m-3-2)

(m-1-2) edge node[right=8pt, below=-2pt] {$ p_1 $} (m-2-2)
(m-2-2) edge node[right=8pt, below=-2pt] {$ p_0 $} (m-3-2)
(m-3-2) edge node[right=8pt, below=-2pt] {$  $} (m-4-2)

(m-2-2) edge node[above] {$ k_{1} $} (m-2-3)
(m-3-2) edge node[above] {$ k_{0} $} (m-3-3)
(m-4-2) edge node[above] {$ $} (m-4-3);
\path[dotted, -stealth]
(m-1-2) edge node[above] {$  \exists ? \  k_{2} $} (m-1-3);
\path[dashed, -stealth]
(m-1-1) edge[bend right = 15] node[below=5pt, left =20pt] {$ \Pnn  $} (m-1-3)
(m-2-1) edge[bend right = 15] node[below=5pt, left =20pt] {$ \Pnn  $} (m-2-3)
(m-3-1) edge[bend right = 15] node[below=5pt, left =20pt] {$ \Pnn  $} (m-3-3);
\end{tikzpicture} 
\end{center}
Similarly, the $k$-invariant $k_2$ exists if and only if the composite $\Pnn k_1$ is nullhomotopic, and the set of such extensions is given by a subset of homotopy classes of maps in $[ \kww1, \Sigma^{3\rn - (2,0)} \HH]$. As for the previous case, the existence of a unique $k$-invariant $k_2$ will be shown in Proposition \ref{prop:existkwn}.

More generally, if the $k$-invariants $k_0, \ldots, k_{m}$ exist, one can define $\kwwmplus$ and ask if we can continue building the tower, i.e., if the next $k$-invariant $k_{m+1}$ exists. The pattern is clear and a $k$-invariant $k_{m+1}$ exists if and only if we have the relation
\begin{equation} \label{eq:existencekinv}
\Pnn k_{m} = 0  \in \HH^{(m+2)\rn-(m,0)}(\kwwm).
\end{equation}
Once $k_{m+1}$ exists, one can alter it by the set of maps
\begin{equation} \label{eq:modulikinv}
\left[ \kwwm, \Sigma^{(m+2)\rn-(m+1,0)} \HH \right] \cong \HH^{(m+2)\rn-(m+1,0)}(\kwwm).
\end{equation}
We will show in Proposition \ref{prop:existkwn} that at each step there exists a unique $k$-invariant $k_m$, which in addition turns out to satisfy the relation $\Pnn k_m = 0$, producing $k_{m+1}$. Having all these $k$-invariants, we can define a motivic spectrum $\kwn$ as the homotopy limit of the tower
\begin{equation} \label{eq:defkwn}
\kwn \coloneqq \holim \left( \ldots \stackrel{p_2}{\lto} \kww1 \stackrel{p_1}{\lto} \kww0  \stackrel{p_0}{\lto} \ast \right).
\end{equation}
The map $\kwn \lto \kww0 = \HH$ will represent the element $1$ in cohomology so we call it the \emph{fundamental class}. Since the choice of $k$-invariants $k_m$ is canonical, this shows that there is a unique homotopy type $\kwn$ whose Postnikov tower has layers $\HH$ that are successively attached via the cohomology operation $\Pnn$. \qed
\end{cons}

Construction \ref{cons:constkwn} contains the framework for the construction of the motivic spectrum $\kwn$. However, as mentioned above, the formula \eqref{eq:defkwn} does not make sense until we show the existence of the $k$-invariants $k_1, k_2, \ldots$, which we do in Proposition \ref{prop:existkwn}. It is folklore that the existence of these $k$-invariants is equivalent to showing that some specific Toda brackets between Eilenberg-Maclane spectra contain the element zero. Although we don't pursue this direction further, the following remark is meant to explain this folklore result.

\begin{remark}[Existence of $k$-invariants from Toda brackets]
The first $k$-invariant $k_1$ exists if and only if the relation $\Pnn k_0=0$ holds in the $\HH$-cohomology of $\kww0 = \HH$. The second $k$-invariant $k_2$ exists if and only if $\Pnn k_1 = 0$ in the cohomology of the 2-stage motivic spectrum $\kww1$. The diagram 
\begin{center}
\begin{tikzpicture}
\matrix (m) [matrix of math nodes, row sep=3em, column sep=5em]
{  \kww0 = \HH & \Sigma^r \HH      & \Sigma^{2r} \HH & \Sigma^{3r} \HH \\
	    & \Sigma^{1,0} \kww1 &                 &  \\
	    & \Sigma^{1,0} \kww0 &                 &  \\
	    & \Sigma^{r + (1,0)} \HH &            &  \\	    	};
\path[thick, -stealth]
(m-1-1) edge node[above] {$ k_0 = \Pnn $} (m-1-2)
(m-1-2) edge node[above] {$ \Pnn $} (m-1-3)
(m-1-3) edge node[above] {$ \Pnn $} (m-1-4)

(m-1-2) edge node[left] {$ i_0 $} (m-2-2)
(m-2-2) edge node[left] {$ p_0 $} (m-3-2)
(m-3-2) edge node[left] {$ k_0 = \Pnn $} (m-4-2);
\path[dashed, -stealth]
(m-2-2) edge node[right = 10pt, below] {$ \exists \ k_{1} $} (m-1-3)
(m-3-2) edge node[right = 50pt, below] {$ \exists \ \langle \Pnn, \Pnn, \Pnn \rangle  $} (m-1-4);
\end{tikzpicture}
\end{center}
shows that the relation $\Pnn k_1 = 0$ holds if and only if the Toda bracket $\langle \Pnn, \Pnn, \Pnn \rangle$ contains an element that is $\Pnn$-divisible. Since the indeterminacy is also $\Pnn$-divisible, this is equivalent to the bracket $\langle \Pnn, \Pnn, \Pnn \rangle$ containing zero. Moreover, the indeterminacy of the Toda bracket corresponds to the choices of such extensions $k_2$. This generalizes to higher Toda brackets for the higher $k_i$'s, but we do not explore this direction, as we will show by other means that $\Pnn k_m = 0$.
\end{remark}

In the following Proposition \ref{prop:existkwn} we now show the existence of these $k$-invariants, as well as their uniqueness up to homotopy.

\begin{prop} \label{prop:existkwn}
There exist unique $k$-invariants $k_0,k_1, k_2, \ldots$ as described in Construction \ref{cons:constkwn}, defining a unique homotopy type $\kwn$ by equation \eqref{eq:defkwn}.
\begin{proof}
Suppose that the tower has been constructed until the stage
\begin{center}
\begin{tikzpicture}
\matrix (m) [matrix of math nodes, row sep=3em, column sep=5em]
{  \Sigma^{(m+1)\rn - (m+1,0)} \HH & \kwwmplus & \Sigma^{(m+2)\rn - (m+1,0)} \HH \\
	\Sigma^{m\rn - (m,0)} \HH & \kwwm & \Sigma^{(m+1)\rn - (m,0)} \HH, \\
    & \vdots &  \\};
\path[thick, -stealth]
(m-1-1) edge node[above] {$ i_m $} (m-1-2)
(m-2-1) edge node[above] {$ i_{m-1} $} (m-2-2)

(m-1-2) edge node[right=14pt, below=-2pt] {$ p_{m} $} (m-2-2)
(m-2-2) edge node[right=18pt, below=-2pt] {$ p_{m-1} $} (m-3-2)

(m-2-2) edge node[above] {$ k_{m} $} (m-2-3);
\path[dotted, -stealth]
(m-1-2) edge node[above] {$  \exists ? \  k_{m+1} $} (m-1-3);
\path[dashed, -stealth]
(m-1-1) edge[bend right = 13] node[below=5pt, left =20pt] {$ \Pnn  $} (m-1-3)
(m-2-1) edge[bend right = 13] node[below=5pt, left =20pt] {$ \Pnn  $} (m-2-3);
\end{tikzpicture} 
\end{center}
and we want to show that there is a unique possible $k$-invariant $k_{m+1}$. Note that in Construction \ref{cons:constkwn} we uniquely constructed the tower in the case $m=0$ so we can start the inductive process. 

We will show in Lemma \ref{lem:Lem1} below that if $\kwwm$ exists (as assumed by the induction hypothesis), then its cohomology $\HHco(\kwwm)$ vanishes in Chow degrees less than  $-m$. This implies that if a $k$-invariant $k_{m+1}$ exists then it is unique, since the set of choices from equation \eqref{eq:modulikinv} is a subset of the cohomology of $\kwwm$ that is concentrated in Chow degrees less than $-m-1$, which vanishes since $-m-1 < -m$.

To show existence, by equation \eqref{eq:existencekinv} we have to show that $\Pnn k_m = 0$. We will show in Lemma \ref{lem:Lem2} by induction that despite the group $\HH^{(m+2)\rn-(m,0)}(\kwwm)$ not vanishing, we still have the relation $\Pnn k_{m} = 0$. By induction, this concludes the existence of unique $k$-invariants $k_0, k_1, \ldots$ and allows us to define a unique homotopy type $\kwn$ by equation \eqref{eq:defkwn}.
\end{proof}
\end{prop}

\begin{lemma} \label{lem:Lem1}
Fix an $m$, and suppose that $\kww0, \kww1, \ldots, \kwwm$ are constructed as in Construction \ref{cons:constkwn}. Then $\HHco(\kwwm)$ vanishes in Chow degrees less than $-m$.
\begin{proof}
We show it by induction. For the initial case recall that $\kww0 = \HH$ and thus $\HHco(\kww0) \cong \Aaa$ is the Steenrod algebra of operations of $\HH$. Recall from Remark \ref{rem:ChowSteerodalg} that it is concentrated in positive Chow degrees, showing the initial case.

Suppose that the Lemma is shown for $\kww0, \ldots, \kwwss$. The cofiber sequence
\begin{equation*}
\Sigma^{s\rn -(s,0)} \HH \lto \kwws \lto \kwwss
\end{equation*}
induces a long exact sequence in cohomology
\begin{equation*}
\cdots \ltoback \Aaa^{(a,b) - s\rn + (s,0)} \ltoback \HH^{a,b}\kwws \ltoback \HH^{a,b} \kwwss \ltoback \cdots.
\end{equation*}
Fix a couple $(a,b)$ in Chow degree less than $-s$, i.e., such that $a - 2b < -s$. Then $\Aaa^{(a,b) - s\rn + (s,0)} =0$ by the initial case since it is in negative Chow degree, and  $\HH^{a,b} \kwwss = 0$ by the inductive hypothesis since it is in Chow degree less than $-s$, which is less than $-(s-1)$. This implies that $\HH^{a,b}\kwws = 0$, showing the inductive step and finishing the proof.
\end{proof}
\end{lemma}

\begin{remark}
This bound is sharp for every $m$, as it can be seen by the long exact sequences in cohomology that the groups $\HHco(\kwwm)$ do not vanish in Chow degree $-m$. The composite $\Pnn k_m \in \HH^{(m+2)\rn-(m,0)}(\kwwm)$ lives in this non-zero group in Chow degree $-m$, so we cannot use Lemma \ref{lem:Lem1} to show that the product $\Pnn k_m$ is zero.
\end{remark}

\begin{lemma} \label{lem:Lem2}
Fix an $m$, and suppose that $\kww0, \kww1, \ldots, \kwwm$ are constructed as in Construction \ref{cons:constkwn}. Then we have the relation \begin{equation*}
\Pnn k_{m} = 0  \in \HH^{(m+2)\rn-(m,0)}(\kwwm).
\end{equation*}
\begin{proof}
Consider the cofiber sequence 
\begin{equation*}
\Sigma^{m\rn - (m,0)} \HH \lto \kwwm \lto \kwwmminus.
\end{equation*}
The natural transformation $\HHco \stackrel{\cdot \Pnn}{\lto} \HH^{(\ast,\ast) + \rn}$ induces a map of long exact sequences, which in bidegree $(\ast,\ast) = (m+1)\rn - (m,0)$ becomes
\begin{center}
\begin{tikzpicture}
\matrix (m) [matrix of math nodes, row sep=3em, column sep=2em]
{  \cdots & \Aaa^{\rn} & \HH^{(m+1)\rn - (m,0)}(\kwwm) & \HH^{(m+1)\rn - (m,0)}(\kwwmminus) & \cdots \\
   \cdots & \Aaa^{2\rn} & \HH^{(m+2)\rn - (m,0)}(\kwwm) & \HH^{(m+2)\rn - (m,0)}(\kwwmminus) & \cdots. \\};
\path[thick, -stealth]
(m-1-2) edge (m-1-1.10)
(m-1-3) edge (m-1-2)
(m-1-4) edge (m-1-3)
(m-1-5) edge (m-1-4)

(m-2-2) edge (m-2-1.10)
(m-2-3) edge (m-2-2)
(m-2-4) edge (m-2-3)
(m-2-5) edge (m-2-4)

(m-1-2) edge node[right] {$ \Pnn \cdot $} (m-2-2)
(m-1-3) edge node[right] {$ \Pnn \cdot $} (m-2-3)
(m-1-4) edge node[right] {$ \Pnn \cdot $} (m-2-4);
\end{tikzpicture} 
\end{center}
Both cohomology groups of $\kwwmminus$ on the right of the above diagram are concentrated in Chow degree $-m$, so they vanish by Lemma \ref{lem:Lem1}. We thus get the commutative square
\begin{center}
\begin{tikzpicture}
\matrix (m) [matrix of math nodes, row sep=3em, column sep=2em]
{  \Pnn  \in \Aaa^{\rn} & \HH^{(m+1)\rn - (m,0)}(\kwwm) \ni k_m \\
   \Pnn\Pnn \in \Aaa^{2\rn} & \HH^{(m+2)\rn - (m,0)}(\kwwm)\ni \Pnn k_m, \\};
\path[thick, -stealth]
(m-1-1) edge node[left] {$ \Pnn \cdot $} (m-2-1)
(m-1-2) edge node[left] {$ \Pnn \cdot $} (m-2-2);
\path[left hook->]
(m-1-2) edge (m-1-1)
(m-2-2) edge (m-2-1);
\end{tikzpicture} 
\end{center}
with injective horizontal maps. The element $k_m$ lives in the top right corner of the square, and by definition is sent to $\Pnn \in \Aaa^{\rn}$ along the top horizontal map. Since $\Pnn\Pnn = 0 \in \Aaa^{2\rn}$ by Lemma \ref{lem:Pstext}, the product $\Pnn k_m$ is also sent to zero in $\Aaa^{2\rn}$ by the bottom horizontal map. Since the bottom horizontal map is injective, it follows that $\Pnn k_m = 0$.
\end{proof}
\end{lemma}

This finishes the construction started in \ref{cons:constkwn}, constructing a homotopy type $\kwn$ defined by its Postnikov tower. From this tower, we will now compute its cohomology. 

\begin{prop} \label{prop:cohomkwn}
The cohomology of $\kwn$ is given as an $\Aaa$-module by
\begin{equation} \label{eq:cohomkwn}
\HHco(\kwn) \cong \Aaa // E(\Pnn).
\end{equation}
In particular, it is concentrated in non-negative Chow degrees.
\begin{proof}
Applying the contravariant functor $\HHco$ to the tower defining $\kwn$ gives an unrolled exact couple and thus an associated spectral sequence as shown in Figure \ref{fig:AHsskwn}.
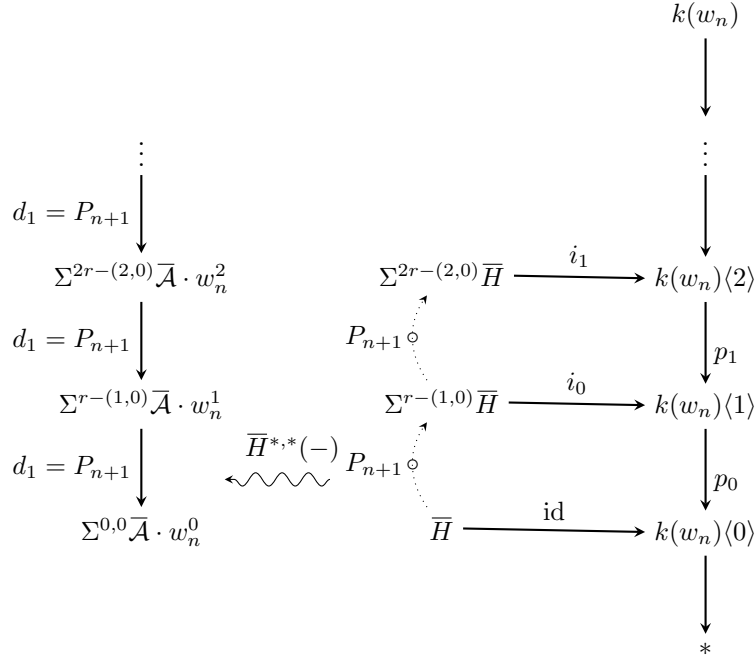
\begin{figure}
\begin{tikzpicture}
\matrix (m) [matrix of math nodes, row sep=3em, column sep=5em]
{ & &  \kwn \\
  \vdots &  & \vdots  \\
  \Sigma^{2r - (2,0)} \Aaa \cdot {w_n^2}	& \Sigma^{2r - (2,0)} \HH & \kww2  \\
	\Sigma^{r - (1,0)} \Aaa \cdot {w_n^1}	& \Sigma^{r - (1,0)} \HH & \kww1 \\
   \Sigma^{0,0} \Aaa \cdot {w_n^0}	& \HH & \kww0 \\
       &  & \ast \\};
\path[thick, -stealth]
(m-3-2) edge node[above] {$ i_1 $} (m-3-3)
(m-4-2) edge node[above] {$ i_0 $} (m-4-3)
(m-5-2) edge node[above] {$ \id $} (m-5-3)

(m-1-3) edge (m-2-3)
(m-2-3) edge (m-3-3)

(m-3-3) edge node[right=8pt, below=-2pt] {$ p_1 $} (m-4-3)
(m-4-3) edge node[right=8pt, below=-2pt] {$ p_0 $} (m-5-3)
(m-5-3) edge node[right=8pt, below=-2pt] {$ $} (m-6-3)

(m-4-1) edge node[left] {$ d_1 = \Pnn  $}  (m-5-1)
(m-3-1) edge node[left] {$ d_1 = \Pnn  $}  (m-4-1)
(m-2-1) edge node[left] {$ d_1 = \Pnn  $}  (m-3-1);

\path [draw,snake it, -stealth]    (-1,-2) -- (-2.4,-2) ;
\node [above] at (-1.5,-1.9) { $\HH^{\ast,\ast}(-)$ };

\path[dotted, -stealth]
(m-4-2) edge[bend left = 35] node[left] {$ \Pnn  $} node  {$ \circ  $} (m-3-2)
(m-5-2) edge[bend left = 35] node[left] {$ \Pnn  $} node  {$ \circ  $}  (m-4-2);
\end{tikzpicture}
\captionof{figure}{The spectral sequence computing the cohomology of $\kwn$ is drawn schematically on the left-hand side.} \label{fig:AHsskwn}
\end{figure}
The $E_{\infty}$-page of this spectral sequence computes the colimit
\begin{equation*}
\mathcal{H}^{\ast,\ast} \coloneqq \colim \left( \HHco(\kww0) \lto \HHco(\kww1) \lto \cdots \right).
\end{equation*}
We will first show that the $E_{\infty}$-page is isomorphic to the quotient $\Aaa//E(\Pnn)$, and then show that in this setting we have
\begin{equation*}
\mathcal{H}^{\ast,\ast} \cong \HH^{\ast,\ast}(\kwn),
\end{equation*}
i.e., that the cohomology of this homotopy limit is computed by the colimit of the cohomologies. We index the $E_1$-page of this spectral sequence as 
\begin{equation*}
E_1^{s,t,w} =  \HH^{t,w} \left( \Sigma^{s\cdot(\rn -(1,0))} \HH \right),
\end{equation*}
where $s$ is the homological degree, and $(t,w)$ are the internal degrees. The first differential $d_1$ is the boundary map in the tower, which is multiplication by the cohomology operation $\Pnn$. Since $\Pnn$ is exact on $\Aaa$ by Corollary \ref{cor:Pstprimextexact}, the $d_1$ differential wipes out everything in homological degree $s >0$, and leaves the quotient $\Aaa//E(\Pnn)$ in degree $s = 0$. It thus follows that the spectral sequence collapses at $E_2$ for degree reasons, with no possible hidden extensions. 

It remains to show that this colimit is isomorphic to $\HHco(\kwn)$. This is shown by a standard technique which we explain for the reader's convenience. Fix a bidegree $(t,w)$ until the end of the proof. Consider the diagram
\begin{center}
\begin{tikzpicture}
\matrix (m) [matrix of math nodes, row sep=3em, column sep=2em]
{  \cdots & \HHtw(\kwwm) & \HHtw(\kwwmplus) & \HHtw(\kwwmplusplus) & \cdots  \\
          & \HHtw(\kwn), & & & \\};
\path[thick, -stealth]
(m-1-1) edge (m-1-2)
(m-1-2) edge (m-1-3)
(m-1-3) edge (m-1-4)
(m-1-4) edge (m-1-5)

(m-1-2) edge (m-2-2)
(m-1-3) edge (m-2-2.35)
(m-1-4.270) edge (m-2-2.05);
\end{tikzpicture} 
\end{center}
where the vertical and diagonal maps are induced by the maps from the diagram defining $\kwn$ as a homotopy limit. For $m$ big enough, the maps in the diagram
\begin{equation*}
\HHtw(\kwwm) \stackrel{\cong}{\lto} \HHtw(\kwwmplus) \stackrel{\cong}{\lto} \cdots \stackrel{\cong}{\lto} \mathcal{H}^{t,w}
\end{equation*}
are all isomorphisms since the successive fibers 
\begin{equation*}
\Sigma^{m\rn - (m,0)} \HH, \ \Sigma^{(m+1)\rn - (m+1,0)} \HH, \ldots
\end{equation*}
are suspended too much and have no cohomology in degree $(t,w)$. Similarly, for $m$ big enough the fiber of $\kwn \lto \kwwm$ is too connected, and thus also has no cohomology. It follows that for $m$ big enough, the maps in the diagram
\begin{center}
\begin{tikzpicture}
\matrix (m) [matrix of math nodes, row sep=3em, column sep=2em]
{  \cdots & \HHtw(\kwwm) & \HHtw(\kwwmplus) & \HHtw(\kwwmplusplus) & \cdots  \\
          & \HHtw(\kwn), & & & \\};
\path[thick, -stealth]
(m-1-1) edge (m-1-2)
(m-1-2) edge node[above] {$ \cong $} (m-1-3)
(m-1-3) edge node[above] {$ \cong $} (m-1-4)
(m-1-4) edge node[above] {$ \cong $} (m-1-5)

(m-1-2) edge node[auto] {$ \cong $} (m-2-2)
(m-1-3) edge node[auto] {$ \cong $} (m-2-2.35)
(m-1-4.270) edge node[auto] {$ \cong $} (m-2-2.05);
\end{tikzpicture} 
\end{center}
become isomorphisms, and thus we have the required isomorphism
\begin{equation*}
\HHtw(\kwn) \stackrel{\cong}{\ltoback} \mathcal{H}^{t,w}.
\end{equation*}
Since $\Aaa$ is concentrated in non-negative Chow degrees by Remark \ref{rem:ChowSteerodalg}, then so is $\HHco(\kwn) \cong \Aaa//E(\Pnn)$.
\end{proof}
\end{prop}

\begin{cor} \label{cor:cohomsmashpowerofkwn}
The cohomology of the smash product $\kwn \smas \kwn$ is given by
\begin{equation*}
\HHco \left( \kwn \smas \kwn \right) \cong \Aaa//E(\Pnn) \otimes \Aaa//E(\Pnn),
\end{equation*}
with the diagonal $\Aaa$-module structure. More generally, for any $m \geq 1$ we have 
\begin{equation*}
\HHco \big( \kwn^{\smas m} \big) \cong \big( \Aaa//E(\Pnn) \big)^{\otimes m},
\end{equation*}
with the iterated diagonal $\Aaa$-module structure. In particular, for any $m$ the cohomology $\HHco \big( \kwn^{\smas m} \big)$ is concentrated in positive Chow degrees.
\begin{proof}
The first part follows from the K\"unneth isomorphism of Proposition \ref{prop:Kunneth}. An easy induction computes the cohomology of the smash product $\kwn^{\smas m}$. Finally, $\Aaa//E(\Pnn)$ is concentrated in positive Chow degrees since $\Aaa$ is (Remark \ref{rem:ChowSteerodalg}), and thus so is $\big( \Aaa//E(\Pnn) \big)^{\otimes m}$.
\end{proof}
\end{cor}

\begin{remark}[Additive description of $\piss(\kwn)$] \label{rem:htpykwn}
It is easy to see from Construction \ref{cons:constkwn} that additively the homotopy groups of $\kwn$ are given by sparse copies of $\F_2$'s
\begin{equation*}
\pisw(\kwn) \cong  \begin{cases} 
   \F_2 & \text{if } (s,w) = m \cdot \rn \text{ for some } m \in \mathbb{N}_0 \\
   0     & \text{otherwise.}
  \end{cases}
\end{equation*} 
For example, one can compute $\piss(\kwn)$ by Milnor's $\lim^1$ short exact sequence, where the $\lim^1$ term vanishes since $p_m$ is surjective in homotopy groups and thus the inverse sequence satisfies the Mittag-Leffler condition. We will compute the ring structure on $\piss(\kwn)$ via the motivic $\HH$-based Adams spectral sequence in Theorem \ref{thm:kwnEoo}.
\end{remark}

\subsection{The $E_{\infty}$ ring structure and $\Kwn$} \label{subsec:Eookwn}

We will now endow the motivic spectrum $\kwn$ with an $E_{\infty}$ ring structure. The technique is the same as in \cite[Section 3]{GheCt} and is done in three steps. The first step is to construct a ring map $\mu \colon \kwn \smas \kwn \lto \kwn$, which we do by lifting the fundamental class $1 \in \HH^{0,0} \big( \kwn^{\smas 2} \big)$ along the Postnikov tower of $\kwn$. We will then show that this endows $\kwn$ with a unital, associative and commutative monoid structure in the homotopy category. The last step is to use Robinson's obstruction theory \cite{RobinsonGamma} to extend it to an $E_{\infty}$ ring structure. 

The following Lemma \ref{lem:smashkwntokwn} provides a homotopy class of maps $\mu \colon \kwn \smas \kwn \lto \kwn$, as well as some estimates necessary to apply Robinson's obstruction theory.

\begin{lemma} \label{lem:smashkwntokwn}
For any $m \geq 1$ the abelian group of homotopy classes of maps satisfies
\begin{equation*}
\left[ \kwn^{\smas m}, \Sigma^{t,w} \kwn \right] \cong  \begin{cases} 
   \F_2 & \text{if } (t,w) = (0,0) \\
   0     & \text{if } (t,w) \text{ is in negative Chow degree, i.e., } t - 2w < 0.
  \end{cases}
\end{equation*}
Moreover, the non-trivial map $\kwn^{\smas m} \lto \kwn$ preserves the fundamental class in $\HH$-cohomology, i.e., sends $1$ to $1^{\otimes m}$.
\begin{proof}
Consider the Atiyah-Hirzebruch spectral sequence computing 
\begin{equation*}
\left[ \kwn^{\smas m}, \kwn \right]_{\ast,\ast}
\end{equation*}
given by applying the functor $\left[ \kwn^{\smas m}, - \right]$ to the tower defining $\kwn$ in Construction \ref{cons:constkwn}. We index the $E_1$-page by
\begin{equation*}
E_1^{s,t,w} = \left[ \kwn^{\smas m}, \Sigma^{t,w} \Sigma^{s\rn - (s,0)} \HH \right] \cong \HH^{(t,w) + s\rn - (s,0)} \big( \kwn^{\smas m} \big),
\end{equation*}
where $s \geq 0$ is the homological degree and $(t,w)$ are the internal degrees. This is a first quadrant spectral sequence (if plotted in the $(s,t)$ plane) and thus converges to
\begin{equation*}
E_{\infty}^{s,t,w} \cong \left[ \kwn^{\smas m}, \Sigma^{t,w} \kwn \right].
\end{equation*}
Since the cohomology of the smash power $\kwn^{\smas m}$ is concentrated in non-negative Chow degrees by Corollary \ref{cor:cohomsmashpowerofkwn}, and since $\rn$ is in Chow degree zero, the $E_1$-page is concentrated in degrees
\begin{equation*}
t - 2w - s \geq 0.
\end{equation*}
In particular, if $t - 2w < 0$, the spot $E_{\infty}^{s,t,w}$ is zero for any $s$ and thus also
\begin{equation*}
\left[ \kwn^{\smas m}, \Sigma^{t,w} \kwn \right] = 0.
\end{equation*}
In the case where $(t,w) = (0,0)$, then necessarily $s=0$. The $l^{\text{th}}$ differential $d_l$ goes from $(s,t,w)$ to $(s+l, t+1,w)$, and there is thus no possible differential entering $E_1^{0,0,0}$. Since $d_l$ reduces the quantity $t - 2w - s$ by $l-1$, the only possible differential exiting $E_1^{0,0,0}$ is a $d_1$. As in Proposition \ref{prop:cohomkwn}, observe that the $d_1$ differential is multiplication by $\Pnn$ on $\HH^{\ast,\ast}\big( \kwn^{\smas m} \big)$. Recall from Corollary \ref{cor:cohomsmashpowerofkwn} that 
\begin{equation*}
\HHco \big( \kwn^{\smas m} \big) \cong \big( \Aaa//E(\Pnn) \big)^{\otimes m},
\end{equation*}
with $\HH^{0,0} \big( \kwn^{\smas m} \big) \cong \Z/2^{\otimes m} \cong \Z/2$. Since $\Pnn \in \Aaa$ is primitive, it acts as zero on $1 \otimes \cdots \otimes 1$. This shows that there is no possible differential on $E_1^{0,0,0}$ and thus
\begin{equation*}
\left[ \kwn^{\smas m}, \kwn \right] \cong E_1^{0,0,0} = \HH^{0,0} \big( \kwn^{\smas m} \big) \cong \Z/2^{\otimes m} \cong \Z/2.
\end{equation*}
Let's call a representative for the non-trivial class by $\mu \colon \kwn^{\smas m} \lto \kwn$. Unwinding this chain of isomorphisms shows that $\mu$  start as the fundamental class $1^{\otimes m} = \mu_0 \colon \kwn^{\smas m} \lto \HH$ and can be uniquely lifted along the Postnikov tower as shown in the diagram
\begin{center}
\begin{tikzpicture}
\matrix (m) [matrix of math nodes, row sep=2em, column sep=5em]
{ &  \kwn & \\
   &  \vdots &  \\
    & \kww2 & \Sigma^{3r - (2,0)} \HH \\
	 & \kww1 & \Sigma^{2r - (1,0)} \HH \\
 \kwn^{\smas m}  & \kww0 & \Sigma^r \HH. \\};
\path[thick, -stealth]
(m-1-2) edge (m-2-2)
(m-2-2) edge (m-3-2)
(m-3-2) edge node[right=8pt, below=-2pt] {$ p_1 $} (m-4-2)
(m-4-2) edge node[right=8pt, below=-2pt] {$ p_0 $} (m-5-2)

(m-3-2) edge node[above] {$ k_{2} $} (m-3-3)
(m-4-2) edge node[above] {$ k_{1} $} (m-4-3)
(m-5-2) edge node[above] {$ k_{0} $} (m-5-3);
\path[dashed, -stealth]
(m-5-1) edge  node[auto] {$ \mu_0 = 1 $} (m-5-2)
(m-5-1) edge  node[auto] {$ \mu_1 $} (m-4-2.190)
(m-5-1) edge[bend left=10] node[auto] {$ \mu_2 $} (m-3-2.190)
(m-5-1) edge[bend left=15] node[auto] {$ \mu $} (m-1-2.200);
\end{tikzpicture}
\end{center}
This shows that $\mu$ sends the fundamental class $1$ to the fundamental class $1^{\otimes m}$ in cohomology since the vertical composite is non-trivial and thus $1$ in cohomology.
\end{proof}
\end{lemma}

\begin{prop} \label{prop:htpymonoid}
The map $\kwn \smas \kwn \stackrel{\mu}{\lto} \kwn$ is homotopy unital, associative and commutative.
\begin{proof}
Since $\pi_{0,0}(\kwn) \cong \F_2$ from Remark \ref{rem:htpykwn}, the non-zero element
\begin{equation*}
\Ss \stackrel{i}{\lto} \kwn
\end{equation*}
will be the unit of the ring structure on $\kwn$. The multiplication $\mu$ is homotopy left unital if and only if the composite
\begin{equation*}
\Ss \smas \kwn \stackrel{i \smas \id}{\ltoo} \kwn \smas \kwn \stackrel{\mu}{\lto} \kwn
\end{equation*}
is homotopy equivalent to the identity map on $\kwn$. By Corollary \ref{cor:cohomsmashpowerofkwn} the group of homotopy classes of self-maps of degree $(0,0)$ on $\kwn$ is $\Z/2$, so it suffices to show that the above composite is not nullhomotopic. We can do so by embedding it in the following commutative diagram
\begin{center}
\begin{tikzpicture}
\matrix (m) [matrix of math nodes, row sep=3em, column sep=3em]
{  \Ss \smas \kwn & \kwn \smas \kwn & \kwn  \\
   \Ss \smas \Ss  & \HH \smas \HH   & \HH, \\};
\path[thick, -stealth]
(m-1-1) edge node[above] {$ i \smas \id $} (m-1-2)
(m-1-2) edge node[above] {$ \mu $} (m-1-3)

(m-2-1) edge node[auto] {$ \id \smas i $} (m-1-1)
(m-1-2) edge node[auto] {$ 1 \smas 1 $} (m-2-2)
(m-1-3) edge node[auto] {$ 1 $} (m-2-3)

(m-2-1) edge node[auto] {$ i \smas i $} (m-2-2)
(m-2-2) edge node[auto] {$ \mu $} (m-2-3);
\end{tikzpicture} 
\end{center}
where both squares are seen to commute by Lemma \ref{lem:smashkwntokwn} since the map $1$ represents 1 in the cohomology of $\kwn$. The top horizontal composite cannot be nullhomotopic since the bottom horizontal composite is the unit of $\HH$ and thus not nullhomotopic. This shows that $\mu$ is left unital, and by a similar argument that $\mu$ is also right unital.

Recall from Corollary \ref{cor:cohomsmashpowerofkwn} that $\mu$ lives in the group $\left[ \kwn \smas \kwn, \kwn \right] \cong \Z/2 \cdot \{ \mu \}$. Precomposing with the factor swap map $\chi \colon \kwn \smas \kwn \lto \kwn \smas \kwn$ is an involution on this group, which forces $\mu \comp \chi \cong \mu$, i.e., showing that $\mu$ is homotopy commutative.

For associativity, we need to compare the two maps $\mu \comp( \mu \smas \id)$ and $\mu \comp ( \id \smas \mu)$ in the group of homotopy classes of maps $\left[ \kwn \smas \kwn \smas \kwn, \kwn \right]$. Since $\mu$ is unital, precomposing both maps with the units
\begin{equation*}
\Ss \smas \Ss \smas \Ss \stackrel{i \smas i \smas i}{\ltooo} \kwn \smas \kwn \smas \kwn
\end{equation*}
gives the non-zero map $\Ss \stackrel{i}{\lto} \kwn$. This means that both maps $\mu \comp( \mu \smas \id)$ and $\mu \comp ( \id \smas \mu)$ are not nullhomotopic, and since $\left[ \kwn \smas \kwn \smas \kwn, \kwn \right] \cong \Z/2$ by Corollary \ref{cor:cohomsmashpowerofkwn}, they are homotopic.
\end{proof}
\end{prop}

\begin{thm} \label{thm:kwnEoo}
For any $n$, the motivic spectrum $\kwn \in \CtCell$ admits an essentially unique $E_{\infty}$ ring structure and satisfies
\begin{itemize}
\item $\HHco(\kwn) \cong \Aaa//E(\Pnn)$ as an $\Aaa$-module,
\item $\piss(\kwn) \cong \F_2[w_n]$ as a ring.
\end{itemize}
\begin{proof}
We start by rigidifying the homotopy ring structure on $\kwn$ to an $E_{\infty}$ ring structure by using Robinson's obstruction theory. This obstruction theory has been adapted to the motivic setting in \cite[Corollary 3.2]{GheCt}. More precisely, the multiplication $\mu$ extends to an $E_{\infty}$ ring structure if the groups
\begin{equation*}
\left[ \kwn^{\smas m}, \Sigma^{3-m',0} \kwn \right]
\end{equation*}
are all zero for $m' \geq 4$ and $2 \leq m \leq m'$. For a fixed $m' \geq 4$, observe that these groups are in negative Chow degree for any $m$, and thus vanish by Lemma \ref{lem:smashkwntokwn}. This shows that the multiplication map $\mu$ can be extended to an $E_{\infty}$ ring structure on $\kwn$. Furthermore, \cite[Corollary 3.2]{GheCt} shows that this $E_{\infty}$ ring structure is unique if 
\begin{equation*}
\left[ \kwn^{\smas m}, \Sigma^{2-m',0} \kwn \right]
\end{equation*}
are all zero for $m' \geq 4$ and $2 \leq m \leq m'$. Another application of Lemma \ref{lem:smashkwntokwn} shows these are zero and thus that $\kwn$ admits a unique $E_{\infty}$ ring structure.

Its cohomology has been computed in Proposition \ref{prop:cohomkwn}. For its homotopy, we already know from remark \ref{rem:htpykwn} that $\piss(\kwn)$ is given by a copy of $\F_2$ in every degree of the form $m \cdot (\rn - (1,0))$ where $m \geq 0$ and $\rn - (1,0)$ is the bidegree of the class $w_n$. To show that the ring structure is polynomial, one can for example consider the $\Ct$-linear $\HH$-based motivic Adams spectral sequence
\begin{equation*}
E_2 = \Ext_{\Aaa}\big( \HHco(\kwn), \F_2 \big) \Longrightarrow \piss(\kwn)
\end{equation*}
which is now multiplicative since $\kwn$ is a motivic ring spectrum. Since $\HHco(\kwn) \cong \Aaa//E(\Pnn)$ we can apply the usual change of rings to the $E_2$-page 
\begin{equation*}
E_2 \cong \Ext_{E(\Pnn)}(\F_2, \F_2) \cong \F_2[w_n].
\end{equation*}
The spectral sequence collapses now at $E_2$ with no possible hidden extensions.
\end{proof}
\end{thm}

Another relation between the homotopy element $w_n$ and the cohomology operation $\Pnn$ is given in the following Proposition.

\begin{prop} \label{prop:cofseqbockstein}
For any $n$, there is a cofiber sequence
\begin{equation*}
\Sigma^{\rn - (1,0)} \kwn \stackrel{w_n}{\lto} \kwn \stackrel{1}{\lto} \HH \stackrel{\beta_n}{\lto} \Sigma^{\rn} \kwn,
\end{equation*}
where the boundary map $\beta_n$ is such that the composite
\begin{equation*}
\HH \stackrel{\beta_n}{\lto} \Sigma^{\rn} \kwn \stackrel{1}{\lto} \Sigma^{\rn} \HH
\end{equation*}
gives the cohomology operation $\Pnn \in \Aaa$.
\begin{proof}
Consider the cofiber sequence
\begin{equation*}
\Sigma^{\rn -(1,0)} \kwn \stackrel{w_n}{\lto} \kwn \lto C \lto \Sigma^{\rn} \kwn,
\end{equation*}
where we denote by $C$ the cofiber of multiplication by $w_n$. Comparing it with the cofiber sequence coming from the beginning of the tower of $\kwn$ gives a diagram
\begin{center}
\begin{tikzpicture}
\matrix (m) [matrix of math nodes, row sep=3em, column sep=3em]
{  \Sigma^{\rn -(1,0)} \kwn & \kwn & C & \Sigma^{\rn} \kwn  \\
   \Sigma^{\rn -(1,0)} \HH & \kww1 & \kww0 = \HH & \Sigma^{\rn} \HH,  \\};
\path[thick, -stealth]
(m-1-1) edge node[above] {$ w_n $} (m-1-2)
(m-1-2) edge node[above] {$ $} (m-1-3)
(m-1-3) edge node[auto] {$ q $} (m-1-4)

(m-2-1) edge node[above] {$ $} (m-2-2)
(m-2-2) edge node[above] {$ 1 $} (m-2-3)
(m-2-3) edge node[auto] {$ \Pnn $} (m-2-4)

(m-1-2) edge node[above] {$ $} (m-2-2)
(m-1-2) edge node[auto] {$ 1 $} (m-2-3);
\end{tikzpicture} 
\end{center}
where $\kwn \lto \kww1$ is the natural map and $1$ denotes the fundamental class. The composite $1 \comp w_n$ lives in the Chow degree $-1$ part of the cohomology of $\kwn$, which is zero by Proposition \ref{prop:cohomkwn}. This implies that there exists a filler $\psi$
\begin{center}
\begin{tikzpicture}
\matrix (m) [matrix of math nodes, row sep=3em, column sep=3em]
{  \Sigma^{\rn -(1,0)} \kwn & \kwn & C & \Sigma^{\rn} \kwn & \Sigma^{1,0} \kwn  \\
   \Sigma^{\rn -(1,0)} \HH & \kww1 & \kww0 = \HH & \Sigma^{\rn} \HH & \Sigma^{1,0} \kww1  \\};
\path[thick, -stealth]
(m-1-1) edge node[above] {$ w_n $} (m-1-2)
(m-1-2) edge node[above] {$ $} (m-1-3)
(m-1-3) edge node[auto] {$ q $} (m-1-4)
(m-1-4) edge node[auto] {$ w_n $} (m-1-5)

(m-2-1) edge node[above] {$ $} (m-2-2)
(m-2-2) edge node[above] {$ 1 $} (m-2-3)
(m-2-3) edge node[auto] {$ \Pnn $} (m-2-4)
(m-2-4) edge node[above] {$ $} (m-2-5)

(m-1-2) edge node[above] {$ $} (m-2-2)
(m-1-2) edge node[auto] {$ 1 $} (m-2-3);
\path[-stealth, dashed]
(m-1-3) edge node[auto] {$ \exists \psi  $} (m-2-3)
(m-1-4) edge node[auto] {$ \exists ! \ \phi  $} (m-2-4);
\end{tikzpicture} 
\end{center}
which itself implies that there is another filler $\phi$ making all squares commute. By the long exact sequence in cohomology, it is easy to see that $\left[ C, \HH \right] \cong \Z/2$. Observe that $\psi$ is non-zero since $\kwn \lto \kww0$ is non-zero, which implies that $\psi$ is unique up to homotopy and induces an isomorphism on homotopy groups. It follows that its cofiber is contractible and thus that it is an equivalence. Observe that $\phi$ is unique up to homotopy since $\left[ \Sigma^{1,0} \kwn, \Sigma^{\rn} \HH \right] = 0$ by Proposition \ref{prop:cohomkwn}. Since $\left[ \Sigma^{1,0} \kwn, \Sigma^{\rn} \HH \right] \cong \Z/2$ and $\Pnn \neq 0$, then $\phi$ is also non-zero, forcing it to be the fundamental class $\phi = 1$. This gives the desired cofiber sequence
\begin{equation*}
\Sigma^{\rn - (1,0)} \kwn \stackrel{ w_n}{\lto} \kwn \stackrel{1}{\lto} \HH \stackrel{\beta_n}{\lto} \Sigma^{\rn} \kwn,
\end{equation*}	
where we denote the composite $q \comp \psi^{-1}$ by $\beta_n$. The composite $\HH \stackrel{\beta_n}{\lto} \Sigma^{\rn} \kwn \stackrel{1}{\lto} \Sigma^{\rn} \HH$ is the cohomology operation $\Pnn$ by the above comparison of cofiber sequences.
\end{proof}
\end{prop}

\begin{cor} \label{cor:Kwn}
For any $n$, there is a motivic $E_{\infty}$ graded field $\Kwn$ with $\piss(\Kwn) \cong \F_2[w_n^{\pm 1}]$.
\begin{proof}
As a module, define $\Kwn$ as the homotopy colimit
\begin{equation*}
\Kwn \coloneqq \hocolim \left( \kwn \stackrel{w_n}{\lto} \Sigma^{-|w_n|} \kwn \stackrel{w_n}{\lto} \cdots \right) .
\end{equation*}
By compactness of $S^{s,w}$, its homotopy groups are given by 
\begin{equation*}
\piss( \Kwn ) \cong \piss( \kwn ) [w_n^{-1}] \cong \F_2[w_n^{\pm 1}].
\end{equation*}
It remains to show that this localization can be performed in motivic $E_{\infty}$ rings. For this, one can apply the methods of \cite{HillHopkins}. We will now give a minimal argument to explain how this applies to the motivic setting, and refer to \cite[Section 3.1]{HillHopkins} for more details. Recall that the motivic $E_{\infty}$ operad that we consider is the simplicial operad where $E \Sigma_n$ is a constant motivic space. In particular, this space admits a cellular filtration where the layers are spheres $S^{\ast,0}$ in weight zero. There is thus a spectral sequence computing the homotopy groups of $(E\Sigma_n)_+ \smas_{\Sigma_n} Z$, which has as input the homotopy groups of various suspensions $\Sigma^{\ast,0} Z$. One can now apply \cite{HillHopkins} since the acyclics form a localizing subcategory that is closed under suspensions of the form $\Sigma^{\ast,0}$.
\end{proof}
\end{cor}


\section{The Motivic Spectrum $wBP$} \label{sec:wBP}

In this section we will construct $E_{\infty}$ ring spectra $\wBP$ and $\wBPn$ with homotopy groups given by
\begin{equation*}
\piss(\wBP) \cong \F_2[w_0, w_1, \ldots] \qquad \text{ and } \qquad \piss(\wBPn) \cong \F_2[w_0, w_1, \ldots, w_n].
\end{equation*}
These spectra will be constructed for the property that
\begin{equation*}
\HHco(\wBP) \cong \Aaa//\EPinf \qquad \text{ and } \qquad \HHco(\wBPn) \cong \Aaa//E(P_1, \ldots, \Pnn)
\end{equation*}
with the natural $\Aaa$-module structure. As we did in Section \ref{subsec:Pst} for $\kwn$, in Section \ref{subsec:formulasforwBP} we derive some formulas in the Steenrod algebra $\Aaa$ and its dual $\Aaadual$. In Section \ref{subsec:wBP} we proceed to construct $\wBP$, by using a version of Toda's Realization Theorem \cite[Lemma 3.1]{Toda}. We finally endow it with an $E_{\infty}$ ring structure in Section \ref{subsec:EoowBP}.

\subsection{More formulas in the $\HH$-Steenrod algebra} \label{subsec:formulasforwBP}

Recall from Proposition \ref{prop:HHSteenrod} that the $\Ct$-linear dual $\HH$-Steenrod algebra, i.e., the Hopf algebra of $\Ct$-linear co-operations on $\HH$ is given by
\begin{equation*}
\piss(\HH \smas_{\Ct} \HH) = \Aaadual \cong \F_2[\xi_1, \xi_2, \ldots] \otimes E(\tau_0, \tau_1, \ldots).
\end{equation*}
In this section we will need to work in Milnor's basis, whose notation we recall.

\begin{nota}[Milnor's basis and the $\PR$ notation] \label{nota:notationPR}
Given a sequence $R=(r_1, r_2, \ldots)$ of non-negative integers with only finitely many non-zero entries, denote by $\PR \in \Aaa$ the dual element to  $\xi_1^{r_1} \xi_2^{r_2} \cdots$. The \emph{length} of a sequence $R$ is the non-negative number
\begin{equation*}
l(R) = r_1 + r_2 + \cdots.
\end{equation*}
Denote by $\Delta_j$ the sequence of length 1 containing a 1 in position $j$, and thus we recover $\P^{\Delta_j} = \Pj$. Given two sequences $R$ and $R'$, denote by 
\begin{equation*}
\P^{R -R'} = \begin{cases} 
   \text{dual to } \xi_1^{r_1 - r_1'} \xi_2^{r_2 -r_2'} \cdots & \text{if } r_j \geq r_j' \text{ for all } j \\
   0      & \text{if not.}
  \end{cases}
\end{equation*}
A sequence $R$ is called \emph{even} if every $r_i$ is even. Given a sequence $R$, denote by $\PRR$ the dual to $\xi_1^{2r_1} \xi_2^{2r_2} \cdots$, and thus $\P^{2\Delta_j}$ is dual to $\xi_j^2$.
\end{nota}

Recall the ungraded injective map of Hopf algebras $\Acl \inj \Aaa$ from equation \eqref{eq:mapofSteenrod} (or graded in the sense of Remark \ref{rem:gradedbyweight}). By using the same $\PR$ notation in the classical setting (and so in $\Acl$ we have $\Pj = Q_{j-1}$), this map sends $\PR$ to $\PR$ and $\Pj = Q_{j-1}$ to $\Pj$. Moreover, the classical formula $c(Q_j) = Q_j \in \Acl$ implies that motivically $c(\Pj) = \Pj \in \Aaa$. 

We will construct $\wBP$ by assembling $\kw0$'s, and since its cohomology is given by $\HHco(\kw0) \cong \Aaa//E(\Pone)$, we need to derive some formulas in the Hopf algebra quotient $\Aaa//E(\Pone)$.

\begin{lemma} \label{lemma:formulasinAmodmodPone}
In $\Aaa//E(\Pone)$, the following relations hold
\begin{enumerate}
\item $\Pjplus = \Pone \cdot \cPdeltaj$, for any $j \geq 2$.
\item $\Pone \cdot \cPR = \displaystyle\sum_{j \geq 1} \cPRdeltaj \cdot \Pone \cdot \cPdeltaj$, for any sequence $R$.
\end{enumerate}
\begin{proof}
In \cite[Section 2]{BP}, it is shown that the following formula 
\begin{equation*}
\PRR Q_0 + Q_0 \PRR = \sum_{j \geq 1} Q_j \PRdeltaj \in \Acl
\end{equation*}
holds in the classical Steenrod algebra for any sequence $R$. We warn the reader that there is a switch in notation between this formula and \cite[Formula 2.5]{BP}, as Brown and Peterson adopt a different notation in the case $p=2$, where they let $\PR$ be the dual of $\xi_1^{2r_1} \xi_2^{2r_2}\cdots$. Through the map $\Acl \inj \Aaa$, this relation gives the motivic formula
\begin{equation} \label{eq:PRbasicformula}
\Pone \PRR + \PRR \Pone = \sum_{j \geq 1} \Pjplus \PRdeltaj.
\end{equation}
When $R = \Delta_j$, this formula becomes
\begin{equation} \label{eq:Pjplus}
\Pone \Pdeltaj + \Pdeltaj \Pone = \Pjplus.
\end{equation}
Applying the anti-morphism $c(-)$ and considering it in the quotient $\Aaa//E(\Pone) = \Aaa \otimes_{E(\Pone)} \F_2$ gives the desired first formula $\Pjplus = \Pone \cPdeltaj$. By plugging equation \eqref{eq:Pjplus} in \eqref{eq:PRbasicformula} we get
\begin{equation*}
\Pone \PRR + \PRR \Pone = \sum_{j \geq 1} \Pone\Pdeltaj\PRdeltaj + \Pdeltaj \Pone \PRdeltaj.
\end{equation*}
Applying the anti-morphism $c(-)$ and considering it in the quotient $\Aaa//E(\Pone) = \Aaa \otimes_{E(\Pone)} \F_2$ gives the desired second formula
\begin{equation*}
\Pone \cdot \cPR = \sum_{j \geq 1} \cPRdeltaj \cdot \Pone \cdot \cPdeltaj.
\end{equation*}
\end{proof}
\end{lemma}

\begin{remark}
The formulas of Lemma \ref{lemma:formulasinAmodmodPone} live in the quotient $\Aaa//E(P_1) \cong \Aaa \otimes_{E(P_1)} \F_2$, where the action of $P_1$ on the \emph{right} is reduced to zero. The first formula \eqref{eq:Pjplus} is symmetric, and the conjugation morphism $c(-)$ is unnecessary for this formula alone. However, the point of applying the conjugation $c(-)$ is because the element $P_{j+1}$ is multiplied on the left in equation \eqref{eq:PRbasicformula}. Without applying the anti-morphism $c(-)$ there would be no simplification after plugging-in \eqref{eq:PRbasicformula} in \eqref{eq:Pjplus} and the formulas would not be as nice. In fact, as we will see in Proposition \ref{prop:layers}, topologically realizing the differential of the chain complex \eqref{eq:algres} would not be possible without applying $c(-)$.
\end{remark}

\subsection{The construction of $\wBP$} \label{subsec:wBP}

We will construct $\wBP$ via a certain tower, in the category of $\Ct$-modules. Recall that we are trying to construct a motivic ring spectrum $\wBP$ whose homotopy groups are given by the polynomial ring
\begin{equation*}
\piss(\wBP) \cong \F_2[w_0, w_1, \ldots],
\end{equation*}
where $w_n$ is an element detected by the cohomology operation $\Pnn$. Unlike the previous section, the Postnikov tower approach is not tractable for $\wBP$ as it is hard to isolate the monomials $w_0^{n_0} w_1^{n_1} \cdots$ of a given bidegree, and thus hard to describe the layers. 

Observe however that as in the case of $\kwn$, it suffices to construct a motivic spectrum with cohomology given by the Hopf algebra quotient $\Aaa//\EPinf$. In fact, by a change of rings theorem and a careful analysis of degrees, the $\Ct$-linear $\HH$-based Adams spectral sequence computing the homotopy of such a motivic spectrum collapses at $E_2 \cong \F_2[w_0, w_1, \ldots]$.

We will construct such a motivic spectrum by following an idea of Toda from \cite[Lemma 3.1]{Toda}, where Toda constructs classical spectra with given cohomology by attaching copies of $H\F_2$ together. We also need to adapt and incorporate a trick which was used in \cite{BP} to construct the classical Brown-Peterson spectrum $\BP$. This trick is to construct $\wBP$ by attaching together wedges of $\kw0$'s, instead of wedges of $\HH$'s. This has the effect of reducing the number of wedge summands in the layers, and also of reducing the complexity of some computations, since these are included in the construction of $\kw0$ that was already done in Theorem \ref{thm:kwnEoo}.

More precisely, we will construct an inverse tower of motivic spectra, whose associated graded will be a topological realization of a resolution of $\Aaa//\EPinf$. The following Proposition is the first step in doing so, by constructing this associated graded, i.e., the layers of the desired tower.

\begin{prop} \label{prop:layers}
There exists a complex of motivic spectra $\kV0 \stackrel{\delta_0}{\lto} \kV1 \stackrel{\delta_1}{\lto} \cdots$, where 
\begin{enumerate}
\item the composite of two consecutive maps is nullhomotopic,
\item each $X_i$ is a locally finite wedge\footnote{i.e., a possibly infinite wedge $\coprod_{\alpha \in A} \Sigma^{r_{\alpha}} X_{\alpha}$ with a finite number of wedge summands $X_{\alpha}$'s in any given (bi)-degree $\rn$.} of suspensions of $\kw0$,
\item the $\HH$-cohomology of this cochain complex is an $\Aaa//E(\Pone)$-free resolution of $\Aaa//\EPinf$.
\end{enumerate}
\begin{proof}

For every $j \geq 2$, consider the periodic bigraded $E(\Pj)$-free resolution of $\Z/2$ given by
\begin{equation} \label{eq:resofexterior}
0 \ltoback \Z/2 \ltoback E(\Pj) \stackrel{\Pj}{\ltoback} \Sigma^{|\Pj|} E(\Pj) \stackrel{\Pj}{\ltoback} \Sigma^{2|\Pj|} E(\Pj) \ltoback \cdots.
\end{equation}
For simplicity, we will denote the exterior subalgebra generated by $\Ptwo, \Pthree, \ldots$ by 
\begin{equation*}
E \coloneqq \EPinf//E(\Pone) \cong E(P_2, P_3, \ldots).
\end{equation*}
By tensoring together these resolutions\footnote{where the term $\Z/2$ is not part of the resolution and $E(\Pj)$ is in homological degree 0.} for every $j \geq 2$, we get a bigraded $E$-free resolution
\begin{equation} \label{eq:algresolutionoverE}
0 \ltoback \Z/2 \ltoback E \otimes V_0 \stackrel{d_0}{\ltoback} E \otimes V_1 \stackrel{d_1}{\ltoback} E \otimes V_2 \ltoback \cdots,
\end{equation}
for some bigraded $\F_2$-vector space $V_i$. A preferred $\F_2$-basis of $V_i$ is given by the set of sequences
\begin{equation} \label{eq:canbasisofVj}
\{ e_R  \  | \ R =(r_2, r_3, \ldots) \text{ satisfies } l(R) = r_2 + r_3 + \cdots = i \}
\end{equation}
of length $i$. 
It follows that the bigrading on $V_i$ is given by
\begin{equation} \label{eq:degreebasisVi}
|e_R| = \sum_{j \geq 2} r_j \cdot | \Pj| = \sum_{j \geq 2} r_j \cdot |\xi_j| = \sum_{j \geq 2} r_j \cdot (2^{j+1} -2, 2^j -1).
\end{equation}
In the notation $E \otimes V_i$, the $E$-linear differential is given on this basis by
\begin{equation*}
d_i(1 \otimes e_R) = \sum_{j \geq 2} \Pj \otimes (e_{R - \Delta_j}).
\end{equation*}
We can now tensor up the resolution of equation \eqref{eq:algresolutionoverE} via $\Aaa//E(\Pone) \otimes_{E} -$ to obtain the algebraic $\Aaa//E(\Pone)$-free resolution 
\begin{equation} \label{eq:algres}
0 \ltoback \Aaa//\EPinf \ltoback \Aaa//E(\Pone) \otimes V_0 \stackrel{d_0}{\ltoback} \Aaa//E(\Pone) \otimes V_1 \stackrel{d_1}{\ltoback} \Aaa//E(\Pone) \otimes V_2 \ltoback \cdots
\end{equation}
of $\Aaa//\EPinf$. The goal is to now realize this resolution topologically. Since the $V_i$ are finite dimensional and the terms in \eqref{eq:algres} are free $\Aaa//E(\Pone)$-modules, they are realized by locally finite wedges of suspensions of $\kw0$ indexed over the same basis. For simplicity, denote by $\kVi$ the bigraded wedge of suspensions of  $k(w_0)$ indexed by the chosen basis of $V_i$ given by equation \eqref{eq:canbasisofVj}, which thus has the prescribed cohomology $\HHco(\kVi) \cong \Aaa//E(\Pone) \otimes V_i$. Similarly, denote by $\HVi$ the bigraded wedge of suspensions of  $\HH$ indexed by the same basis of $V_i$. To realize the differentials, observe that by Lemma \ref{lemma:formulasinAmodmodPone}, the differential $d_i$ in \eqref{eq:algres} can be simplified to the formula
\begin{equation} \label{eq:formuladiff}
d_i(1 \otimes e_R) = \sum_{j \geq 2} \Pj \otimes (e_{R - \Delta_j}) = \Pone \cdot \sum_{j \geq 2} \cPdeltajminus \otimes (e_{R - \Delta_j}).
\end{equation}
The differential $d_i$ can thus be realized by the composite
\begin{equation*}
\delta_i \colon \kVi \stackrel{1}{\lto} \HVi \stackrel{\coprod \cPdeltaj}{\ltoooo} \Sigma^{-(2,1)} \HViplus \stackrel{\beta_0}{\lto} \kViplus,
\end{equation*}
where $\beta_0$ is the Bockstein from Proposition \ref{prop:cofseqbockstein}. The middle map is a locally finite matrix with entries in $\Aaa$, where for a given sequence $R$ of length $i$, it is assembled from the maps
\begin{equation*}
\cPdeltaj \colon \HH\{e_{R}\} \lto \HH\{e_{R + \Delta_j}\} \inj \HViplus.
\end{equation*}
The composite $\delta_i$ realizes the differential $d_i$ since for a given sequence $R$ of length $i+1$ we have
\begin{center}
\begin{tikzpicture}
\matrix (m) [matrix of math nodes, row sep=4em, column sep=6em]
{  \kVi & \HVi & \Sigma^{-(2,1)} \HViplus & \kViplus  \\
        &      &                          & \HH\{e_{R}\}, \\};
\path[thick, -stealth]
(m-1-1) edge node[above] {$ 1 $} (m-1-2)
(m-1-2) edge node[above] {$ \coprod \cPdeltaj $} (m-1-3)
(m-1-3) edge node[auto] {$ \beta_0 $} (m-1-4)

(m-1-4) edge node[auto] {$ 1 $} (m-2-4);
\path[-stealth, dashed]
(m-1-3) edge node[auto] {$ \Pone $} (m-2-4);
\end{tikzpicture} 
\end{center}
which recovers exactly formula \eqref{eq:formuladiff} since $1 \comp \beta_0 = \Pone$ by Proposition \ref{prop:cofseqbockstein}. We have thus defined a sequence of motivic spectra
\begin{equation*}
\kV0 \stackrel{\delta_0}{\lto} \kV1 \stackrel{\delta_1}{\lto} \kV2 \stackrel{\delta_2}{\lto} \cdots,
\end{equation*}
in which each term is a locally finite wedge of suspensions of $\kw0$'s, and which produces an $\Aaa//E(\Pone)$-free resolution of $\Aaa//\EPinf$ after applying $\HH$-cohomology.

It remains to show that the composites $\delta_{i+1} \comp \delta_i$ are nullhomotopic. This is accomplished by the following commutative diagram
\begin{center}
\begin{tikzpicture}
\matrix (m) [matrix of math nodes, row sep=5em, column sep=3.5em]
{  \kVi & \HVi & \Sigma^{-(2,1)} \HViplus & \kViplus  & \HViplus & \HViplusplus  & \kViplusplus, \\};
\path[thick, -stealth]
(m-1-1) edge node[above] {$ 1 $} (m-1-2)
(m-1-2) edge node[above] {$ \coprod \cPdeltaj $} (m-1-3)
(m-1-3) edge node[auto] {$ \beta_0 $} (m-1-4)

(m-1-4) edge node[auto] {$ 1 $} (m-1-5)
(m-1-5) edge node[auto] {$ \coprod \cPdeltaj $} (m-1-6)
(m-1-6) edge node[auto] {$ \beta_0 $} (m-1-7);
\path[-stealth, dashed]
(m-1-1.50) edge[bend left=30] node[auto] {$ \delta_i $} (m-1-4)
(m-1-4.50) edge[bend left=30] node[auto] {$ \delta_{i+1} $} (m-1-7)
(m-1-3) edge[bend right=25] node[below] {$ \Pone $} (m-1-5)
(m-1-1) edge[bend right=30] node[below] {$ \Pone \comp \coprod \cPR \comp 1 $} (m-1-6);
\end{tikzpicture} 
\end{center}
where the wedges $\coprod \cPR$ are taken over sequences $R$ of length 2, and the composite $\kVi \lto \HViplusplus$ is identified from a sum over all such $R$'s of the equation
\begin{equation*}
\Pone \cPR = \sum_{j \geq 1} \cPRdeltaj \Pone \cPdeltaj \in \Aaa//E(\Pone)
\end{equation*}
from lemma \ref{lemma:formulasinAmodmodPone}. The total composite 
\begin{equation*}
\beta_0 \comp \Pone \comp \coprod \cPR \comp 1 = 0
\end{equation*}
is zero since $\beta_0 \comp \Pone = \beta_0 \comp 1 \comp \beta_0$, and $\beta_0 \comp 1 = 0$ since they are consecutive maps in the cofiber sequence of Proposition \ref{prop:cofseqbockstein}.
\end{proof}
\end{prop}

The next step is to construct an inverse tower of motivic spectra, whose layers and induced $d_1$-differential are exactly the cochain complex of motivic spectra
\begin{equation*}
\kV0 \stackrel{\delta_0}{\lto} \kV1 \stackrel{\delta_1}{\lto} \kV2 \stackrel{\delta_2}{\lto} \cdots
\end{equation*}
from Proposition \ref{prop:layers}. The idea is that once we construct this tower, we can compute the cohomology of its inverse limit by the spectral sequence emerging from applying the cohomological functor $\HHco(-)$. We will define $\wBP$ to be the inverse limit of this tower. The $E_1$-page of the associated spectral sequence is the cohomology of the layers, i.e., the cohomology of the above cochain complex. We just showed in Proposition \ref{prop:layers} that this cohomology forms a resolution of $\Aaa//\EPinf$, and thus the spectral sequence collapses at $E_2 = E_{\infty}$ with output $\Aaa//\EPinf$.

\begin{cons}[Construction of $\wBP$ via Toda's realization method] \label{cons:wBP}
Recall that the goal is now to construct a tower of motivic spectra 
\begin{equation*}
X_{-1} = \ast \ltoback X_0 \ltoback X_1 \ltoback \cdots,
\end{equation*}
with layers and induced $d_1$-differential given by the cochain complex
\begin{equation*}
\kV0 \stackrel{\delta_0}{\lto} \kV1 \stackrel{\delta_1}{\lto} \kV2 \stackrel{\delta_2}{\lto} \cdots.
\end{equation*}
The beginning of the tower is given by 
\begin{center}
\begin{tikzpicture}
\matrix (m) [matrix of math nodes, row sep=4em, column sep=5em]
{         		& \kV0 		   \\
  X_{-1} = \ast & X_0 = \kV0                \\
 \Sigma^{1,0} \kV0 &  \kV1,   \\};
\path[thick, -stealth]
(m-2-1) edge node[auto] {$ k_{-1} $} (m-3-1)
(m-2-2) edge node[above] {$ p_{-1} $} (m-2-1)
(m-1-2) edge node[auto] {$ i_{-1} = \id $} (m-2-2);
\path[-stealth, dashed]
(m-1-2) edge[bend right=50] node[left=4pt, above =20pt] {$ \delta_0 $} (m-3-2);
\path[-stealth, dotted]
(m-2-2) edge node[auto] {$ \exists ? \ k_{0} $} (m-3-2);
\end{tikzpicture} 
\end{center}
where the suspension by $\Sigma^{1,0}$ of the first layer is a small adjustment to get the correct output. Evidently, since $i_{-1} = \id$, there is a unique filler $k_0$ up to homotopy. We can thus set $k_0 = \delta_0$, denote its fiber by $X_1$, and ask if the following filler $k_1$ exists in the diagram
\begin{center}
\begin{tikzpicture}
\matrix (m) [matrix of math nodes, row sep=4em, column sep=5em]
{         		& \kV0 		 & \Sigma^{-1,0} \kV1  \\
  X_{-1} = \ast & X_0 = \kV0 &  X_1                \\
 \Sigma^{1,0} \kV0 &  \kV1   & \Sigma^{-1,0} \kV2.  \\};
\path[thick, -stealth]
(m-2-1) edge node[auto] {$ k_{-1} $} (m-3-1)
(m-2-2) edge node[auto] {$ k_{0} = \delta_0 $} (m-3-2)

(m-2-2) edge node[above] {$ p_{-1} $} (m-2-1)
(m-2-3) edge node[above] {$ p_{0} $} (m-2-2)

(m-1-2) edge node[auto] {$ i_{-1} = \id $} (m-2-2)
(m-1-3) edge node[auto] {$ i_{0} $} (m-2-3);
\path[-stealth, dashed]
(m-1-2) edge[bend right=50] node[left=4pt, above =20pt] {$ \delta_0 $} (m-3-2)
(m-1-3) edge[bend right=50] node[left=5pt, above =20pt] {$ \delta_1 $} (m-3-3);
\path[-stealth, dotted]
(m-2-3) edge node[auto] {$ \exists ? \ k_1 $} (m-3-3);
\end{tikzpicture} 
\end{center}
By taking the fiber of $i_0$, this problem becomes an extension problem in the cofiber sequence
\begin{center}
\begin{tikzpicture}
\matrix (m) [matrix of math nodes, row sep=3em, column sep=4em]
{  \Sigma^{-1,0} X_0 & \Sigma^{-1,0} \kV1 & X_1 &  X_0  \\
                      &  &\Sigma^{-1,0} \kV2, & \\};
\path[thick, -stealth]
(m-1-1) edge node[above] {$ k_0 = \delta_0 $} (m-1-2)
(m-1-2) edge node[above] {$ i_0 $} (m-1-3)
(m-1-3) edge node[above] {$ p_0 $} (m-1-4)
(m-1-2) edge node[below=8pt, left] {$ \delta_1  $} (m-2-3);
\path[dotted, -stealth]
(m-1-3) edge node[right] {$  \exists ? \  k_{1} $} (m-2-3);
\end{tikzpicture}
\end{center}
where a filler $k_1$ exists since the composite $\delta_1 \delta_0$ is nullhomotopic by Propositon \ref{prop:layers}. The choices of such extensions are parametrized by the quotient
\begin{equation*}
\quotient{ \left[ X_0, \Sigma^{-1,0} \kV2 \right] }{ \delta_0\text{-divisible elements.}}
\end{equation*}
Recall that both $X_0 = \kV0$ and $\kV2$ are locally finite wedges of suspension of $\kw0$'s, and we can read from equation \eqref{eq:degreebasisVi} that the bidegrees of all suspensions are in Chow degree 0. The set of homotopy classes of maps $\left[ X_0, \Sigma^{-1,0} \kV2 \right]$ is thus built out of self-maps of $\kw0$ in Chow degree $-1$. This is zero by Lemma \ref{lem:smashkwntokwn}, and thus there exists a unique filler $k_1$. 

Getting this far was the base case for the inductive process. Suppose now that the tower
\begin{center}
\begin{tikzpicture}
\matrix (m) [matrix of math nodes, row sep=4em, column sep=4em]
{       &  	      & \Sigma^{-(n-1),0} \kVnmoins & \Sigma^{-n,0} \kVn      \\
 \cdots & X_{n-2} & X_{n-1}                     & X_n  \\
        & \Sigma^{-(n-2),0} \kVnmoins & \Sigma^{-(n-1),0} \kVn & \Sigma^{-n,0} \kVnplus, \\};
\path[thick, -stealth]
(m-2-2) edge node[auto] {$ k_{n-2} $} (m-3-2)
(m-2-3) edge node[auto] {$ k_{n-1} $} (m-3-3)

(m-2-2) edge (m-2-1)
(m-2-3) edge node[above] {$ p_{n-2} $} (m-2-2)
(m-2-4) edge node[above] {$ p_{n-1} $} (m-2-3)

(m-1-3) edge node[auto] {$ i_{n-2} $} (m-2-3)
(m-1-4) edge node[auto] {$ i_{n-1} $} (m-2-4);
\path[-stealth, dashed]
(m-1-3) edge[bend right=40] node[left=8pt, above =20pt] {$ \delta_{n-1} $} (m-3-3)
(m-1-4) edge[bend right=40] node[left=5pt, above =20pt] {$ \delta_{n} $} (m-3-4);
\path[-stealth, dotted]
(m-2-4) edge node[auto] {$ \exists ? \ k_n $} (m-3-4);
\end{tikzpicture} 
\end{center}
has been constructed for some $n \geq 2$. As above, we can desuspend one step, and rewrite this extension problem as
\begin{center}
\begin{tikzpicture}
\matrix (m) [matrix of math nodes, row sep=3em, column sep=4em]
{  \Sigma^{-1,0} X_{n-1} & \Sigma^{-n,0} \kVn & X_n &  X_{n-1}  \\
                      &  &\Sigma^{-n,0} \kVnplus. & \\};
\path[thick, -stealth]
(m-1-1) edge node[above] {$ k_{n-1} $} (m-1-2)
(m-1-2) edge node[above] {$ i_{n-1} $} (m-1-3)
(m-1-3) edge node[above] {$ p_{n-1} $} (m-1-4)
(m-1-2) edge node[below=8pt, left] {$ \delta_n  $} (m-2-3);
\path[dotted, -stealth]
(m-1-3) edge node[right] {$  \exists ? \  k_{n} $} (m-2-3);
\end{tikzpicture}
\end{center}
We will now show both that $\delta_n k_{n-1}=0$ and thus that a filler $k_n$ exists; and that $\left[ X_{n-1}, \Sigma^{-n,0} \kVnplus \right] = 0$ and thus that such a filler is unique. 

Let's first deal with the uniqueness part, as the statement we show will come up in the existence part as well. A slightly more general result is true, namely that the $\kw0$-cohomology of any $X_m$ vanishes in Chow degrees strictly smaller than $-m$. This is easy to show by an induction on $m$ (for all $m$ less than $n$, so that $X_m$ is already constructed), and is completely analogous to Lemma \ref{lem:Lem1}. The base case is $X_0 = \kV0$, whose $\kw0$-cohomology vanishes in negative Chow degrees by Lemma \ref{lem:smashkwntokwn}. The induction is done by inspecting the long exact sequence in cohomology of the cofiber sequence $\Sigma^{-m,0} \kVm \lto X_m \lto X_{m-1}$. We refer to Lemma \ref{lem:Lem1} for more details. We can now use this statement to show that $\left[ X_{n-1}, \Sigma^{-n,0} \kVnplus \right] = 0$. In fact, this set of maps is made out of the part in Chow degree $-n$ of the $\kw0$-cohomology of $X_{n-1}$, which vanishes since $-n < -(n-1)$.

For the existence part, consider the $\Sigma^{-1,0}$-desuspension of the cofiber sequence
\begin{equation*}
\Sigma^{-(n-1),0} \kVnmoins \stackrel{i_{n-2}}{\lto} X_{n-1} \stackrel{p_{n-2}}{\lto}  X_{n-2}.
\end{equation*}
By applying the cohomological functor $\left[ -, \Sigma^{-n,0} \kVnplus \right]$, we get a long exact sequence
\begin{equation*}
\cdots \ltoback \left[ \Sigma^{-n,0} \kVnmoins, \Sigma^{-n,0} \kVnplus \right] \stackrel{i_{n-2}^{\ast}}{\ltoback} \left[ \Sigma^{-1,0} X_{n-1}, \Sigma^{-n,0} \kVnplus \right] \stackrel{p_{n-2}^{\ast}}{\ltoback}  \left[ \Sigma^{-1,0} X_{n-2}, \Sigma^{-n,0} \kVnplus \right] \ltoback \cdots,
\end{equation*}
which contains the composite $\delta_n k_{n-1}$ in its middle term. The right term $\left[ X_{n-2}, \Sigma^{-(n-1),0} \kVnplus \right]$ vanishes since it is concentrated in the Chow degree $-(n-1)$ part of the $\kw0$-cohomology of $X_{n-2}$. This simplifies the long exact sequence to 
\begin{equation*}
\left[ \Sigma^{-n,0} \kVnmoins, \Sigma^{-n,0} \kVnplus \right] \stackrel{i_{n-2}^{\ast}}{\ltoback} \left[ \Sigma^{-1,0} X_{n-1}, \Sigma^{-n,0} \kVnplus \right] \ltoback 0.
\end{equation*}
The image of $\delta_n k_{n-1}$ under the precomposition map $i_{n-2}^{\ast}$ is
\begin{equation*}
i_{n-2}^{\ast}(\delta_n k_{n-1}) = \delta_n k_{n-1} i_{n-2} = \delta_n \delta_{n-1},
\end{equation*}
which is zero by Proposition \ref{prop:layers}. Since $i_{n-2}^{\ast}$ is injective, it follows that $\delta_n k_{n-1} = 0$ and thus that $k_n$ exists. \qed
\end{cons}

\begin{defn}
Define a motivic spectrum $\wBP \in \CtCell$ as the inverse limit of the tower
\begin{equation} \label{eq:defwBP}
\wBP \coloneqq \holim \big( \ast = X_{-1} \ltoback X_0 \ltoback X_1 \ltoback \cdots \big)
\end{equation}
from Construction \ref{cons:wBP}.
\end{defn}

It remains to show that this motivic spectrum has the correct cohomology. The situation is very similar to the case of $\kwn$ in Section \ref{sec:kwn}.

\begin{prop} \label{prop:cohomwBP}
The cohomology of $\wBP$ is given as an $\Aaa$-module by the Hopf algebra quotient
\begin{equation*}
\HHco \left( \wBP \right) \cong \Aaa//\EPinf.
\end{equation*}
\begin{proof}
This is analogous to Proposition \ref{prop:cohomkwn}, which we refer to for more details. Applying the contravariant functor $\HHco$ to the tower defining $\wBP$ gives a spectral sequence computing the algebraic colimit $\colim \HHco(X_i)$. The $E_1$-page is given by the cohomology of the layers $\kVi$, which form a resolution of $\Aaa//\EPinf$ by Proposition \ref{prop:layers}. The $E_2$-page is thus given by $\Aaa//\EPinf$ concentrated in homological degree $s=0$, and the spectral sequence collapses with no possible hidden extensions. Observe that the layers $\Sigma^{-(i-1),0} \kVi$ are more and more connected. In fact, the element of lowest bidegree corresponds to the sequence $R=(r_2=i, 0,0, \ldots)$, and so 
$\Sigma^{-(i-1),0} \kVi$ has no cohomology in degrees lower than $(5i+1,3i)$. By using this bound, the exact same proof as in Proposition \ref{prop:cohomkwn} shows that
\begin{equation*}
\HHco(\wBP) \stackrel{\cong}{\ltoback} \colim \HHco(X_i) \cong \Aaa//\EPinf.
\end{equation*}
\end{proof}
\end{prop}

\begin{cor}
For any $m \geq 1$ we have an isomorphism of $\Aaa$-modules
\begin{equation*}
\HHco \big( \wBP^{\smas m} \big) \cong \big( \Aaa//\EPinf \big)^{\otimes m},
\end{equation*}
where the right hand side has the iterated diagonal $\Aaa$-module structure. In particular, for any $m$ the cohomology $\HHco \big( \wBP^{\smas m} \big)$ is concentrated in positive Chow degrees.
\begin{proof}
This proof is similar to the proof of Corollary \ref{cor:cohomsmashpowerofkwn}.
\end{proof}
\end{cor}


\subsection{The $E_{\infty}$ ring structure on $\wBP$, and $\wBPn$} \label{subsec:EoowBP}

Having constructed a motivic spectrum $\wBP \in \CtCell$ with correct cohomology, we will now follow the same methodology as we did for $\kwn$ in Section \ref{subsec:Eookwn}. We first construct a ring map $\mu \colon \wBP^{\smas 2} \lto \wBP$ by understanding some parts of various groups of homotopy classes of maps $\left[ \wBP^{\smas m}, \wBP \right]_{\ast,\ast}$. We then show that $\mu$ turns $\wBP$ into a unital, associative and commutative monoid in the homotopy category. Finally, the vanishing of $\left[ \wBP^{\smas m}, \wBP \right]_{\ast,\ast}$ in some particular degrees feeds Robinson's obstruction theory \cite{RobinsonGamma} which rigidifies $\mu$ to an $E_{\infty}$ ring structure on $\wBP$. The $\Ct$-linear $\HH$-based motivic Adams spectral sequence is now multiplicative and collapses at $E_2$, showing that $\wBP$ has polynomial homotopy in the periodicity elements $w_0, w_1, \ldots$.

Most of these steps are very similar to the case of $\kwn$. We will only sketch the arguments in these cases and refer to the appropriate proof in Section \ref{sec:kwn} for more details. Since $\wBP$ is built out of a tower whose layers are wedges of suspensions of $\kw0$'s, there is a spectral sequence computing $\left[ \wBP^{\smas m}, \wBP \right]_{\ast,\ast}$ whose layers are made out of $\left[ \wBP^{\smas m}, \kw0 \right]_{\ast,\ast}$. We first say something about these layers.

\begin{lemma} \label{lemma:wBPsmasninkw0}
For any $m \geq 1$ we have
\begin{equation*}
\left[ \wBP^{\smas m}, \Sigma^{t,w} \kwn \right] \cong  \begin{cases} 
   \F_2 & \text{if } (t,w) = (0,0) \\
   0     & \text{if } (t,w) \text{ is in negative Chow degree, i.e., } t - 2w < 0.
  \end{cases}
\end{equation*}
\begin{proof}
The exact same proof as in Lemma \ref{lem:smashkwntokwn} applies by changing $\kw0$ to $\wBP$, since the cohomology of the smash powers $\wBP^{\smas m}$ is also concentrated in non-negative Chow degrees, with a copy of $\Z/2$ in degree $(0,0)$.
\end{proof}
\end{lemma}

\begin{lemma} \label{lemma:wBPsmasninwBP}
For any $m \geq 1$ we have
\begin{equation*}
\left[ \wBP^{\smas m}, \Sigma^{t,w} \wBP \right] \cong  \begin{cases} 
   \F_2 & \text{if } (t,w) = (0,0) \\
   0     & \text{if } (t,w) \text{ is in negative Chow degree, i.e., } t - 2w < 0.
  \end{cases}
\end{equation*}
Moreover, the non-trivial map $\wBP^{\smas m} \lto \wBP	$ sends $1$ to $1^{\otimes m}$ in cohomology.
\begin{proof}
As in Lemma \ref{lem:smashkwntokwn} and Lemma \ref{lemma:wBPsmasninkw0}, this is another Atiyah-Hirzebruch spectral sequence, by applying the functor $\left[ \wBP^{\smas m} , - \right]$ to the inverse limit tower defining $\wBP$. With the indexing
\begin{equation*}
E_1^{s,t,w} = \left[ \wBP^{\smas m}, \Sigma^{t,w} \Sigma^{-s,0} \kVs \right] \cong \left[ \wBP^{\smas m}, \Sigma^{t-s,w} \kVs \right],
\end{equation*}
this is a first quadrant spectral sequence in the $(s,t)$ plane which converges to 
\begin{equation*}
E_{\infty}^{s,t,w} \cong \left[ \wBP^{\smas m}, \Sigma^{t,w} \wBP \right].
\end{equation*}
Recall that each $\kVs$ is a locally finite wedge of suspensions of $\kw0$, where the suspensions are in Chow degree 0. By Lemma \ref{lemma:wBPsmasninkw0}, the $\kw0$-cohomology of $\wBP$ is concentrated in non-negative Chow degree, with only a $\Z/2$ in degree $(0,0)$. The proof of Lemma \ref{lem:smashkwntokwn} applies by changing $\HH$ to $\kw0$, and $\kw0$ to $\wBP$.
\end{proof}
\end{lemma}

In the case $m=2$, call a representative of the non-trivial class of maps by $\mu \colon \wBP \smas \wBP \lto \wBP$. We are now ready to show that this map can be extended to an $E_{\infty}$ ring structure on $\wBP$.

\begin{thm} \label{thm:wBP}
The motivic spectrum $\wBP \in \CtCell$ admits an essentially unique $E_{\infty}$ ring structure and satisfies
\begin{itemize}
\item $\HHco(\wBP) \cong \Aaa//\EPinf$ as an $\Aaa$-module,
\item $\piss(\wBP) \cong \F_2[w_0,w_1,\ldots]$ as a ring.
\end{itemize}
\begin{proof}
First use the proof of Proposition \ref{prop:htpymonoid} by changing $\kw0$ to $\wBP$ to show that $\mu$ turns $\wBP$ into a unital, associative and commutative monoid in the homotopy category. Similarly, by using the vanishing results of Lemma \ref{lemma:wBPsmasninwBP}, the proof of Theorem \ref{thm:kwnEoo} shows that $\mu$ can be uniquely extended to an $E_{\infty}$ ring structure. Finally, the $\Ct$-linear $\HH$-based motivic Adams spectral sequence 
\begin{equation*}
E_2 = \Ext_{\Aaa}\big( \HHco(\wBP), \F_2 \big) \Longrightarrow \piss(\wBP)
\end{equation*}
is multiplicative since $\wBP$ is a ring spectrum. Since $\HHco(\wBP) \cong \Aaa//\EPinf$ we can apply the usual change of rings to the $E_2$-page 
\begin{equation*}
E_2 \cong \Ext_{\EPinf}(\F_2, \F_2) \cong \F_2[w_0,w_1,\ldots].
\end{equation*}
The spectral sequence collapses now at $E_2$ with no possible hidden extensions.
\end{proof}
\end{thm}

\begin{cor} \label{cor:wBPn}
For any $n \geq 0$ there exists a motivic spectrum $\wBPn \in \CtCell$, which admits an essentially unique $E_{\infty}$ ring structure and satisfies
\begin{itemize}
\item $\HHco(\wBPn) \cong \Aaa//E(\Pone, \Ptwo, \ldots, \Pnplus)$ as an $\Aaa$-module,
\item $\piss(\wBPn) \cong \F_2[w_0,w_1,\ldots, w_n]$ as a ring.
\end{itemize}
\begin{proof}
The whole proof is very similar to the case of $\wBP$, so we will only indicate what needs to be changed. Tensoring together for $2 \leq j \leq n+1$ the resolutions of equation \eqref{eq:resofexterior} gives a resolution
\begin{equation*}
0 \ltoback \Z/2 \ltoback E \otimes V_0 \stackrel{d_0}{\ltoback} E \otimes V_1 \stackrel{d_1}{\ltoback} E \otimes V_2 \ltoback \cdots,
\end{equation*}
where $E = E(\Ptwo, \Pthree, \ldots, \Pnplus)$ and where $V_i$ is a finite dimensional bigraded $\F_2$-vector space with basis given by the set of sequences
\begin{equation*} 
\{ R =(r_2, r_3, \ldots, r_{n+1}) \  | \ l(R) = r_2 + r_3 + \cdots + r_{n+1} = i \}.
\end{equation*}
Tensoring up to $\Aaa//E(\Pone) \otimes_{E} -$ gives the $\Aaa//E(\Pone)$-free resolution
\begin{equation*}
0 \ltoback \Aaa//E(\Ptwo, \Pthree, \ldots, \Pnplus) \ltoback \Aaa//E(\Pone) \otimes V_0 \stackrel{d_0}{\ltoback} \Aaa//E(\Pone) \otimes V_1 \stackrel{d_1}{\ltoback} \Aaa//E(\Pone) \otimes V_2 \ltoback \cdots.
\end{equation*}
The differential is still given by 
\begin{equation*}
d_i(1 \otimes R) = \sum_{j \geq 2} \Pj \otimes (R - \Delta_j) = \Pone \cdot \sum_{j \geq 2} \cPdeltajminus \otimes (R - \Delta_j),
\end{equation*}
since the formula we use still lives in the quotient $\Aaa//E(\Pone)$. From here on, the rest of the proof is a copy of the proof for $\wBP$, with the advantage that the terms $\kVi$ are now finite wedges of suspensions of copies of $\kw0$'s.
\end{proof}
\end{cor}

\begin{remark}
One can also construct the underlying motivic spectrum of $\wBPn$ by taking the quotient
\begin{equation*}
\quotient{\wBP}{w_{n+1},w_{n+2}, \ldots}.
\end{equation*}
The fact that these quotients are $E_{\infty}$-rings requires some extra work from this point of view though.
\end{remark}



\bibliographystyle{alpha}
\bibliography{mybibliography}

\end{document}